\documentclass[letter,12pt]{elsarticle}
\usepackage{graphicx}
\usepackage{epsfig}

\usepackage{amsmath,amsfonts,amssymb,amsthm}
\usepackage{latexsym,graphicx}

\theoremstyle{plain}
\newtheorem{thm}{Theorem}
\newtheorem*{thm*}{Theorem}

\newtheorem{lem}{Lemma}
\newtheorem*{lem*}{Lemma}
\newtheorem{cor}{Corollary}
\newtheorem{conj}{Conjecture}
\newtheorem{prop}{Proposition}

\theoremstyle{definition}

\newcommand{\be}{\begin{equation}}
\newcommand{\ee}{\end{equation}}

\newcommand{\1}{{\rm 1\hspace*{-0.4ex}%
\rule{0.1ex}{1.52ex}\hspace*{0.2ex}}}

\def\cR{{\mathcal{R}}}
\def\cS{{\mathcal{S}}}
\def\cE{{\mathcal{E}}}

\journal{Chaos, Solitons \& Fractals}

\begin{document}

\begin{frontmatter}

\title{Tokunaga and Horton self-similarity \\
for level set trees of Markov chains}
\author[IZ]{Ilia Zaliapin\corref{cor}}
\ead{zal@unr.edu}
\cortext[cor]{Corresponding author, Phone: 1-775-784-6077, Fax: 1-775-784-6378}
\author[YK]{Yevgeniy Kovchegov}
\ead{kovchegy@math.oregonstate.edu}
 
\address[IZ]{Department of Mathematics and Statistics,
University of Nevada, Reno, NV, 89557-0084, USA}
\address[YK]{Department of Mathematics, 
Oregon State University,
Corvallis, OR,
97331-4605, USA}

\begin{abstract}
The Horton and Tokunaga branching laws provide a convenient 
framework for studying self-similarity in random trees. 
The Horton self-similarity is a weaker property that 
addresses the {\it principal branching} in a tree; 
it is a counterpart of the power-law size distribution
for elements of a branching system.
The stronger Tokunaga self-similarity 
addresses so-called {\it side branching}.
The Horton and Tokunaga self-similarity have been empirically 
established in numerous observed and modeled systems, 
and proven for two paradigmatic models:
the critical Galton-Watson branching process with finite progeny and
the finite-tree representation of a regular Brownian excursion.
This study establishes the Tokunaga and Horton self-similarity 
for a tree representation of a finite symmetric homogeneous 
Markov chain.
We also extend the concept of Horton and Tokunaga self-similarity
to infinite trees and establish self-similarity for an infinite-tree
representation of a regular Brownian motion.
We conjecture that fractional Brownian motions are also Tokunaga
and Horton 
self-similar, with self-similarity parameters depending on
the Hurst exponent. 
\end{abstract}


\end{frontmatter}

\section{Introduction and motivation}
\label{intro}

Hierarchical branching organization is ubiquitous in nature. 
It is readily seen in river basins, drainage networks, bronchial 
passages, botanical trees, and snowflakes, to mention but a few 
(e.g., \cite{Shreve66,NTG97,TPN98,NBW06}). 
Empirical evidence reveals a surprising similarity among various 
natural hierarchies --- many of them are closely approximated by 
so-called {\it self-similar trees} (SSTs)
\cite{Shreve66,NTG97,TPN98,BWW00,Horton45,Strahler,
Shreve69,Tok78,TBR88,VG00,DR00,PT00,Oss92,Pec95,PG99}.
An SST preserves its statistical structure, in a sense to be defined, 
under the operation of {\it pruning}, i.e., cutting the leaves; this is why the SSTs
are sometimes referred to as {\it fractal} trees \cite{NTG97}.
A two-parametric subclass of {\it Tokunaga} SSTs, 
introduced by Tokunaga \cite{Tok78} in a hydrological context, 
plays a special role in theory and applications, 
as it has been shown to emerge in unprecedented 
variety of modeled and natural phenomena.
The Tokunaga SSTs with a broad range of parameters are seen 
in studies of river networks \cite{Shreve66,BWW00,Shreve69,Tok78,TBR88,Pec95,OWZ97}, 
vein structure of botanical leaves \cite{NTG97,TPN98},
numerical analyses of diffusion limited aggregation \cite{Oss92,MT93},
two dimensional site percolation \cite{TMM+99,YNT+05,ZWG06a,ZWG06b},
and nearest-neighbor clustering in Euclidean spaces \cite{Webb09}.
The diversity of these processes and models hints at the existence 
of a universal (not problem-specific) underlying mechanism 
responsible for the Tokunaga self-similarity and prompts the question: 
{\it What probability models may produce Tokunaga self-similar trees?}
An important answer to this question was given by 
Burd {\it et al.} \cite{BWW00} who studied Galton-Watson 
branching processes and have shown that, in this class, the 
Tokunaga self-similarity is a characteristic property of 
a {\it critical binary branching}, that is the discrete-time 
process that starts with a single progenitor and 
whose members equiprobably either split in two or die at every step.
The critical binary Galton-Watson process is equivalent to 
the Shreve's random river network model, for which the Tokunaga 
self-similarity has been known for long time  
\cite{Shreve66,BWW00,Shreve69,Pec95}.
The Tokunaga self-similarity has also been rigorously 
established in a general hierarchical coagulation 
model of Gabrielov {\it et al.} \cite{GNT99} introduced 
in the framework of self-organized criticality, 
and in a random self-similar network model of Veitzer and 
Gupta \cite{VG00} developed as an alternative to the Shreve's
random network model for river networks.

Prominently, the results of Burd {\it et al.} \cite{BWW00} 
reveal the Tokunaga self-similarity for 
any process represented by the finite Galton-Watson 
critical binary branching. 
In the context of this paper, the most important
example is a regular Brownian motion, whose various connections 
to the Galton-Watson processes are well-known (see Pitman 
\cite{Pitman} for a modern review). 
For instance, the topological structure of the so-called 
$h$-excursions of a regular Brownian motion \cite{NP89} and
a Poisson sampling of a Brownian excursion \cite{Hobson}
are equivalent to a finite critical binary Galton-Watson tree
(Sect.~\ref{TF} below explains the tree representation
of time series), 
and hence these processes are Tokunaga self-similar.

This study further explores Tokunaga self-similarity by focusing 
on trees that describe the topological structure of the level 
sets of a time series or a real function, 
so-called {\it level-set trees}.
Our set-up is closely related to the classical Harris 
correspondence between trees and finite random walks \cite{Harris}, 
and its later ramifications that include infinite trees 
with edge lengths \cite{BWW00,OWZ97,Pitman,Ald1,Ald2,Ald3,LeGall93,LeGall05}.
The main result of this paper is the Tokunaga and closely 
related Horton self-similarity for the level-set trees of 
finite symmetric homogeneous Markov chains (SHMCs) ---
see Sect.~\ref{main}, Theorem~\ref{T1}.
Notably, the Tokunaga and Horton self-similarity concepts have 
been defined so far only for finite trees (e.g., \cite{BWW00,Pec95,MG08}).   
We suggest here a natural extension of Tokunaga and Horton self-similarity
to infinite trees and establish self-similarity for an infinite-tree 
representation of a regular Brownian motion.
The suggested approach is based on the {\it forest of trees 
attached to the floor line} as described by Pitman \cite{Pitman}.
Finally, we discuss the strong distributional self-similarity 
that characterizes Markov chains with exponential jumps.

The paper is organized as follows.
Section~\ref{trees} introduces  
planar rooted trees, trees with edge lengths, Harris paths, 
and spaces of random trees with the Galton-Watson distribution.
The trees on continuous functions are described in Sect.~\ref{TF}. 
Several types of self-similarity for trees --- 
Horton, Tokunaga, and distributional self-similarity --- are 
discussed in Sect.~\ref{SS}.
The main results of the paper are summarized in Sect.~\ref{main}.
Section~\ref{add} addresses special properties of exponential
Markov chains that, in particular, enjoy the strong 
distributional self-similarity.
Proofs are collected in Sect.~\ref{proofs}.
Section~\ref{discussion} concludes.

\section{Trees}
\label{trees}
We introduce here planar trees, the corresponding Harris paths,
and the space of Galton-Watson trees following
Burd {\it et al.} \cite{BWW00}, Ossiander {\it et al.} \cite{OWZ97}
and Pitman \cite{Pitman}. 

\subsection{Planar rooted trees}
Recall that a {\it graph} $\mathcal{G}=(V,E)$ is a collection of vertices
(nodes) 
$V=\{v_i\}$, $1\le i \le N_V$ and edges (links) 
$E=\{e_k\}$, $1\le k \le N_E$. 
In a {\it simple} graph each edge is defined as an unordered
pair of distinct vertices: 
$\forall\, 1\le k \le N_E, \exists! \, 1\le i,j \le N_V, i\ne j$ 
such that $e_k=(v_i,v_j)$ and
we say that the edge $k$ {\it connects} vertices $v_i$ and $v_j$.
Furthermore, each pair of vertices in a simple graph may have 
at most one connecting edge.
  
\begin{figure}[p] 
\centering\includegraphics[width=.5\textwidth]{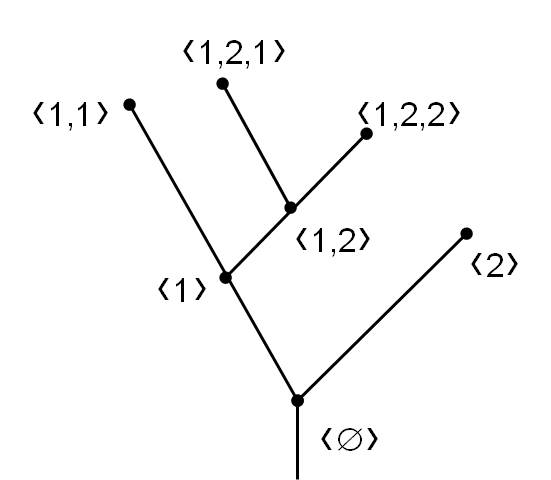}
\caption[Representation of a tree]
{Representation of a tree via a set of finite sequences
$\langle i_1,\dots,i_n\rangle$.}
\label{fig_seq}
\end{figure}

A {\it tree} is a connected simple graph $T=(V,E)$ without cycles,
which readily gives $N_E=N_V-1$.
In a {\it rooted} tree, one node is designated as a root; this
imposes a natural {\it direction} of edges as well as the 
parent-child relationship between the vertices. 
Specifically, we follow \cite{BWW00} to represent a labeled 
(planar) tree $T$ rooted at $\phi$ by a bijection between the 
set of vertices $V$ and set of finite integer-valued sequences
$\langle i_1,\dots,i_n \rangle\in T$ such that
\begin{itemize}
\item[(i)] $\phi=\langle\emptyset\rangle$,
\item[(ii)] if $\langle i_1,\dots,i_n\rangle\in T$ then
$\langle i_1,\dots,i_k\rangle\in T\quad\forall\,1\le k\le n$,
and 
\item[(iii)] if $\langle i_1,\dots,i_n\rangle\in T$ then
$\langle i_1,\dots,i_{n-1}, j\rangle\in T\quad\forall\,1\le j\le i_n$.
\end{itemize}
This representation is illustrated in Fig.~\ref{fig_seq}.
If $v = \langle i_1,\dots,i_n\rangle\in T$ then
$u = \langle i_1,\dots,i_{n-1}\rangle\in T$ is called 
the {\it parent} of $v$, and $v$ is a {\it child} of $u$.
A {\it leaf} is a vertex with no children.
The number of children of a vertex 
$u=\langle i_1,\dots,i_n\rangle\in T$ equals to
$c(u)=\max\{j\}$ over such $j$ that 
$\langle u,j\rangle\equiv\langle i_1,\dots,i_n,j\rangle\in T$.
A {\it binary} labeled rooted tree is represented 
by a set of binary sequences with elements $i_k=1,2$, 
where 1,2 represent the left and right planar directions,
respectively.  
Two trees are called {\it distinct} if they are represented 
by distinct sets of the vertex-sequences.
We complete each tree $T$ by a special {\it ghost edge} $\epsilon$
attached to the root $\phi$, so each vertex in the tree has 
a single parental edge.
A natural direction of edges is from a vertex $v$ to its 
parent $v_p$. 
 
In these settings, the total number of distinct trees with $n$ 
leaves, according to the Cayley's formula, is $n^{n-2}$.
The total number of distinct binary trees with $n$ leaves
is given by the $(n-1)$-th Catalan number \cite{Pitman}
\[C_{n-1}=\frac{1}{n}\left(\begin{array}{c} 2n-2\\ n-1 \end{array}\right).\]

\subsection{Trees with edge-lengths and Harris path}
A tree with {\it edge-lengths} $T=(V,E,W)$ assigns a positive 
lengths $w(e)$ to each edge $e$, $W=\{w(e)\}$; such trees are also
called {\it weighted trees} (e.g., \cite{BWW00,OWZ97}).
The sum of all edge lengths is called the {\it tree length};
we write $\textsc{length}(T)=\sum_e\,w(e)$.  
We call the pair $(V,E)$ a {\it combinatorial tree} and
write $(V,E)=\textsc{shape}(T)$, emphasizing that the
lengths are disregarded in this representation.

\begin{figure}[p] 
\centering\includegraphics[width=0.7\textwidth]{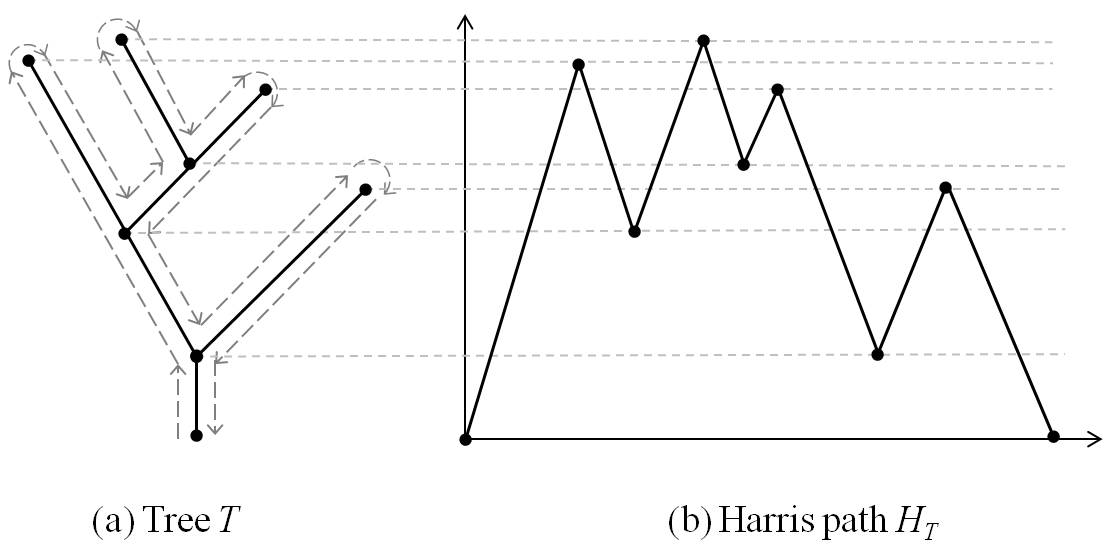}
\caption[Harris path]
{(a) Tree $T$ and its depth-first search illustrated by dashed arrows.
(b) Harris path for the tree $T$ of panel (a). }
\label{fig_Harris}
\end{figure}

If a tree is represented graphically in a plane,
there is a unique continuous map
\[\sigma_T\,:\,[0,\,2\textsc{length}(T)]\to\,T\]
that corresponds to the {\it depth-first search} of $T$,
illustrated in Fig.~\ref{fig_Harris}(a).
The depth-first search starts at the root of planar tree 
with edge-lengths and contours it, moving at a unit speed,  
from left to right so that each edge is traveled twice --- 
its left side in a move away from the root, while its right 
side in a move towards the root.
The {\it Harris path} for a tree $T$ is a continuous function
$H_T(s)\,:\,[0,2\textsc{length}(T)]\to\mathbb{R}$
that equals to the distance from the root 
traveled along the tree $T$ in the depth-fist search.
Accordingly, for a tree $T$ with $n$ leaves, the Harris path 
$H_T(s)$ is a continuous excursion --- 
$H_T(0)=H_T(2\textsc{length}(T))=0$ and $H_T(s)> 0$
for any $s\in \left(0,2\textsc{length}(T)\right)$ ---
that consists of $2n$ linear segments of alternating slopes 
$\pm 1$ \cite{Pitman}, as illustrated in Fig.~\ref{fig_Harris}(b).
The closely related {\it Harris walk} $H_n(k)$, $0\le k \le 2n$ 
for a tree with $n$ vertices is defined as a 
linearly interpolated discrete excursion with $2n$ steps that 
corresponds to the depth-first search that marks each vertex 
in a tree \cite{Harris,Pitman}.
Clearly, the Harris path and Harris walk, as functions
$[0,2\textsc{length}(T)]\to\mathbb{R}$,
have the same trajectory. 
A binary tree with $n$ leaves has $2n-1$ vertices; accordingly,
its Harris path consists of $2n$ segments, and its Harris walk
consists of $4n-2=2(2n-1)$ steps.  

\subsection{Galton-Watson trees}
\label{GW}
The space $\mathbb{T}$ of planar rooted trees with metric
\[d(\tau,\psi)=\frac{1}{1+\sup\{n:\tau|n=\psi|n\}},\]
where $\tau|n=\{\langle i_1,\dots,i_k\rangle \in \tau:k\le n\}$
form a Polish metric space, with the countable dense
subset $\mathbb{T}_0$ of finite trees \cite{OWZ97,BWW00}.  
An important, and most studied, class of distributions 
on $\mathbb{T}$ is the {\it Galton-Watson distribution};
it corresponds to the trees generated by the Galton-Watson
process with a single progenitor and the branching distribution
$\{p_k\}$.
Formally, the distribution $GW_{\{p_k\}}$ assign the following
probability to a closed ball $B\left(\tau,1/n\right)$,
$\tau\in\mathbb{T}$, $n=1,2,\dots$:
\[{\sf P}\left( B\left(\tau,\frac{1}{n}\right)\right)
=\prod_{v\in\tau|(n-1)} p_{c(v)},\]
where $c(v)$ is the number of children of vertex $v$ \cite{BWW00,OWZ97}.  

The classical work of Harris \cite{Harris} notices that the Harris walk 
for a Galton-Watson tree with unit edge-lengths, $n$ vertices and geometric
offspring distribution is an unsigned excursion of length 
$2n$ of a random walk with independent steps $\pm 1$.
Hence, by the conditional Donsker's theorem \cite{Pitman}, 
a properly normalized Harris walk should weakly converge 
to a Brownian excursion.
Aldous \cite{Ald1,Ald2,Ald3}, LeGall \cite{LeGall93,LeGall05}, and
Ossiander {\it et al.} \cite{OWZ97}
have shown that the same limiting behavior is seen
for a broader class of Galton-Watson trees, which may 
have non-trivial edge-lengths and non-geometric offspring 
distribution.
\begin{thm}{\rm \cite[Theorem 3.1]{OWZ97}}
Let $T_n$ be a Galton-Watson tree with the 
total progeny $n$ and offspring distribution
$L$ such that 
gcd$\{j:{\sf P}(L=j)>0\}=1$,
${\sf E}(L)=1$, and
$0<{\sf Var}(L)=\sigma^2<\infty$,
where
gcd$\{\cdot\}$ denotes the greatest common divisor.
Suppose that the i.i.d. lengths $W=\{w(e)\}$ are positive,
independent of $T_n$, have mean 1 and variance $s^2$ and
assume that 
$\lim_{x\to\infty}(x\,\log x)^2{\sf P}(|w(\phi)-1|>x)=0.$
Then the scaled Harris walk $H_n(k)$ converges in distribution
to a standard Brownian excursion $B^{\rm ex}_t$: 
\[\{H_n(2nt)/\sqrt{n},0\le t\le 1\}
\overset{d}{\to} \{2\sigma^{-1}\,B^{\rm ex}_t,0\le t\le 1\},
\quad{\rm as~} n\to\infty.\]
\end{thm}

This paper explores an ``inverse'' problem --- it describes 
trees that correspond to a given finite or infinite Harris 
walk.
We show, in particular, that the class of trees that correspond
to the Harris walks that weakly converge to a Brownian excursion
$B^{\rm ex}_t$ is much broader than the space of Galton-Watson 
trees.

\section{Trees on continuous functions}
\label{TF}
Let $X_t\equiv X(t)\in C\left([L,R]\right)$ be a continuous
function on a finite interval $[L,R]$, $L,R<\infty$.
This section defines the tree associated with $X_t$.
We start with a simple situation when $X_t$ has a 
finite number of local extrema and continue with general 
case. 

\subsection{Tamed functions: Level set trees}
Suppose that the function $X_t\in C\left([L,R]\right)$ has
a finite number of local extrema.
The level set $\mathcal{L}_{\alpha}\left(X_t\right)$ is defined 
as the pre-image of the function values above $\alpha$: 
\[\mathcal{L}_{\alpha}\left(X_t\right) = \{t\,:\,X_t\ge\alpha\}.\]
The level set $\mathcal{L}_{\alpha}$ for each $\alpha$ is
a union of non-overlapping intervals; we write 
$|\mathcal{L}_{\alpha}|$ for their number.
Notice that 
(i) $|\mathcal{L}_{\alpha}| = |\mathcal{L}_{\beta}|$ 
as soon as the interval $[\alpha,\,\beta]$ does not contain a value of
local minima of $X_t$, 
(ii) $|\mathcal{L}_{\alpha}| \ge |\mathcal{L}_{\beta}|$ for any
$\alpha > \beta$, and
(iii) $0\le |\mathcal{L}_{\alpha}| \le n$, where $n$ is the number 
of the local maxima of $X_t$. 

\begin{figure}[p] 
\centering\includegraphics[width=0.5\textwidth]{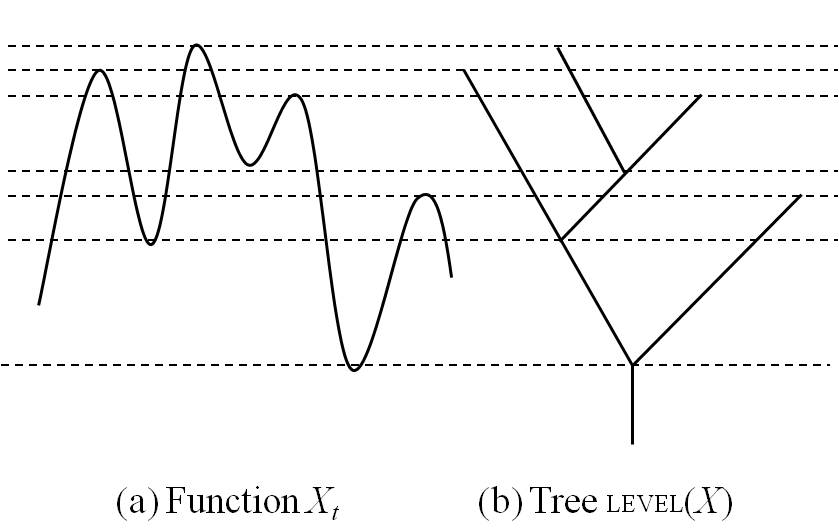}
\caption{Function $X_t$ (panel a) with a finite number of local
extrema and its level-set tree $\textsc{level}(X)$ (panel b).
}
\label{fig3}
\end{figure}

The {\it level set tree} $\textsc{level}(X_t)$ describes the topology 
of the level sets $\mathcal{L}_{\alpha}$ as a function of
threshold $\alpha$, as illustrated in Fig.~\ref{fig3}.
Namely, there are bijections between 
(i) the leaves of $\textsc{level}(X_t)$ and the local
maxima of $X_t$,
(ii) the internal (parental) vertices of $\textsc{level}(X_t)$ 
and the local minima of $X_t$ (excluding possible local minima
at the boundary points), and
(iii) the edges of $\textsc{level}(X_t)$ 
and the first positive excursions of $X(t)-X(t_i)$ to 
right and left of each local minima $t_i$.
The leftmost and rightmost edges $\langle 1,1,\dots,1\rangle$
and $\langle 2,2,\dots,2\rangle$ may correspond to meanders,
that is to a positive segments of $X(t)-X(t_i)$, rather than to 
excursions.
It is readily seen that any function $X_t$ with distinct 
values of the local minima corresponds to a binary tree 
$\textsc{level}(X_t)$.
In this case, the bijection (iii) can be separated into
the bijections between
(iii\,a) the edges $\langle \dots, 1\rangle$ of $\textsc{level}(X_t)$  
and the first positive excursions of $X(t)-X(t_i)$ to the left 
of each local minima $t_i$, and
(iii\,b) the edges $\langle \dots, 2\rangle$ of $\textsc{level}(X_t)$  
and the first positive excursions of $X(t)-X(t_i)$ to the right
of each local minima $t_i$.
The edge $e=(v,u)$ that connects the vertices $v$ and $u$ 
is assigned the length $w(e)$ equal to the absolute difference 
between the values of the respective local extrema of $X_t$
--- according to the bijections (i), (ii) above.

To complete the above construction, a special care should be taken 
of the edge $\epsilon$ attached to the tree root. 
Specifically, let $t_i$, $i=1,\dots,n$, be the set of {\it internal} 
local minima of $X_t$, defined as the set of points such that for any 
$i$ there exists such an open interval $(a_i,b_i)\ni t_i$ that 
$X(t_i)\le X(s)$ for any $s\in(t_i,b_i)$, $X(t_i)<X(b_i)$, and 
$X(t_i)<X(s)$ for any $s\in(a_i,t_i)$.
The last definition treats only the leftmost point of any constant-level
through as a local minima. 
The root of the tree $\textsc{level}(X_t)$ corresponds
to the lowest internal minimum.
If the global minimum $M$ of $X_t$ is reached at one of the 
boundary points, say at $X(L)$, the root of $\textsc{level}(X_t)$ has 
the parental edge $\epsilon$ with the length
$w(\epsilon) = \min_i\left(X(t_i)\right)-X(L)$.
At the same time, if the global minimum $M$ of $X_t$, 
is reached at one of the internal local minima, that is if 
$M = \min_i\left(X(t_i)\right)<\min\left(X(L),X(R)\right)$,
then $|\mathcal{L}_{\alpha}| =0$ for any $\alpha<M$ and
$|\mathcal{L}_{\alpha}| >1$ for any $\alpha>M$.
In other words, the root of  
$\textsc{level}(X_t)$ does not have the parental edge.
In this case, we add the ghost parental edge $\epsilon$ with edge length
$w(\epsilon)=1$.
We write $\textsc{level}(X_t,w(\epsilon))$
to explicitly indicate the length of the ghost edge that might
be added to the level-set tree and save notation $w(\epsilon)$ 
for the value defined above uniquely for each function $X_t$.

By construction, the level 
set trees are invariant with respect to monotone transformations 
of time and values of $X_t$:
\begin{prop}
\label{inv1}
Let $F(\cdot)$ and $G(\cdot)$ be monotone functions
such that 
$Y_t=F\left(X_{G(t)}\right)$ is a continuous function on 
$G\left([L,R]\right).$
Then the function $Y_t$ has the same
combinatorial level set tree as the original 
function $X_t$, that is
\[\textsc{shape}\left(\textsc{level}(X_t,1)\right)=
\textsc{shape}\left(\textsc{level}(Y_t,1)\right).\]
\end{prop}

The tree with edge lengths 
$\textsc{level}(X_t,1)$
is completely specified by the set of the local extrema 
of $X_t$ and its boundary values, and 
is independent of the detailed structure of the intervals 
of monotonicity.
To formalize this observation, 
we write $\cE_X(s)$ for the {\it linear extreme function} obtained 
from $X_t$ by 
(i) linearly interpolating its consecutive local extrema and the 
two boundary values, and 
(ii) changing time within each monotonicity interval as to have
only constant slopes $\pm 1$.
The function $\cE_X(s)$ hence is a piece-wise linear function
with slopes $\pm 1$.
The length of the domain of this function equals the total 
variation of $X_t$.   
We shift this domain to start at 
$s_0 = w(\epsilon) + X(L) - \min_i\left(X(t_i)\right)$,
where $t_i$ are the points of internal local minima as defined 
above.

\begin{prop}
\label{inv2}
The level set tree of a function $X_t$ coincides with 
that of the linear extreme function $\cE_X$:
$\textsc{level}(X_t,1)
=\textsc{level}\left(\cE_X,1\right).$
\end{prop}

The particular domain specification of $\cE_X(z)$ is explained by 
the following statement. 
\begin{prop}
\label{Harris}
Let $H_T(s)$, $s\in[0,2\textsc{length}(T)]$ be the Harris path of the 
level set tree $T=\textsc{level}(X_t,1)$, then
$H_T(z) = \cE_X(z)$ on the domain $D$ of $\cE_X$. 
The domains of $H_T(z)$ and $\cE_X(z)$ coincide, 
i.e. $D=[0,2\textsc{length}(T)]$, if and only if
$X_t$ is a positive excursion, and $D\subset [0,2\textsc{length}(T)]$
otherwise.   
\end{prop}

It is known that each piece-wise linear positive excursion
(Harris path) that consists of $2n$ segments with slopes $\pm 1$ uniquely 
specifies a tree $T$ with no vertices of degree 2 (e.g., \cite{Pitman}). 
Recall that a Harris path corresponds to the depth-first search
that visits each edge in a tree twice; hence the Harris path $H_T$
over-specifies the corresponding tree $T$.
Similarly, the function $\cE_X(s)$ uniquely specifies (and, probably,
over-specifies) the tree
$\textsc{level}(X_t,1)$ with no vertices of degree 2.
If $X_t$ has distinct values of the local minima, then
$\cE_X(s)$ uniquely specifies the binary tree $\textsc{level}(X_t,1)$.

Our definition of the level-set tree cannot be directly
applied to a continuous function with infinite number of
local extrema, say to a trajectory of a Brownian motion.
This motivates the general set-up reviewed in the next 
section \cite{Pitman,LeGall93}.

\subsection{General case}
\begin{figure}[p] 
\centering\includegraphics[width=0.7\textwidth]{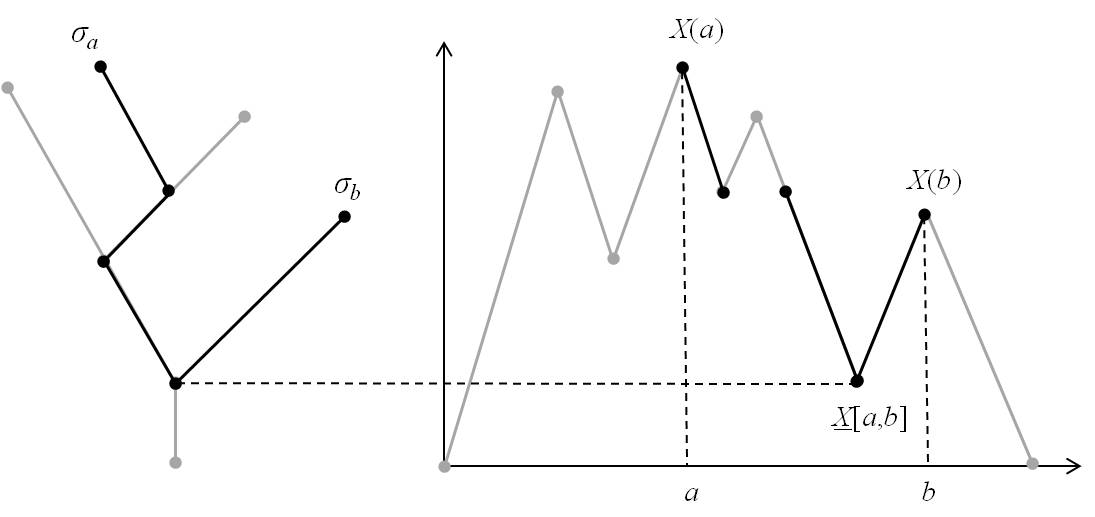}
\caption{Illustration of the pseudo-distance $d_X(a,b)$ used 
to define the tree $\textsc{tree}(X)$ for a continuous 
function $X_t$. 
This example refers to a Harris path $X_t$ 
with the finite number of extrema, so one can 
construct a level set tree for $X_t$. 
Here, the local maxima $X(a)$ and $X(b)$ correspond 
to the leaves $\sigma_a$ and $\sigma_b$ in 
the tree shown on the left. 
The distance between these points is measured along
the shortest path from $\sigma_a$ to $\sigma_b$ along
the tree (marked by heavy lines), or equivalently, 
by Eq. \eqref{tree_dist}.}
\label{fig_tree}
\end{figure}

Let $X_t\equiv X(t)\in C\left([L,R]\right)$ and
$\underline{X}[a,b]:=\inf_{t\in[a,b]} X(t)$, for any
$a,b\in\,[L,R]$.
We define a {\it pseudo-metric} on $[L,R]$ as
\be
\label{tree_dist}
d_X(a,b):=\left(X(a)-\underline{X}[a,b]\right)
+ \left(X(b)-\underline{X}[a,b]\right),\quad a,b\in\,[L,R].
\ee
It is easily verified that if $X_t$ is the Harris
path for a finite tree $T$ and $\sigma_T$ is the corresponding 
depth-first search, then $d_X(a,b)$ equals the distance along 
the tree $T$ between the points $\sigma_T(a)$ and $\sigma_T(b)$
(see Fig.~\ref{fig_tree}).
We write $a\sim_X b$ if $d_X(a,b)=0$.
Accordingly, we define tree $\textsc{tree}(X)$ for the 
function $X_t$ as the metric space $\left([L,R]/\sim_X,d_X\right)$
\cite{Pitman}.  

\vspace{.5cm}
{\bf Remark.}
The definition of the level set tree can be readily applied
to a real-valued Morse function $f:\mathbb{M}\to\mathbb{R}$ 
on a smooth manifold $\mathbb{M}$.
This is convenient for studying functions in higher-dimensional 
domains; see, for instance, Arnold \cite{Arnold} and 
Edelsbrunner {\it et al.} \cite{EHZ03}.
The Harris-path and metric-space definitions are not readily
applicable to multidimensional domains.

\section{Self-similar trees}
\label{SS}
This section describes the three basic forms of the tree 
self-similarity: 
(i) Horton laws, (ii) Self-similarity of side-branching, 
and (iii) Tokunaga self-similarity.
They are based on the Horton-Strahler and Tokunaga schemes
for ordering vertices in a rooted binary tree.
The presented approach was introduced by Horton \cite{Horton45} for 
ordering hierarchically organized river tributaries; 
the methods was later refined by Strahler \cite{Strahler} 
and further expanded by Tokunaga \cite{Tok78} to include so-called
side-branching.

\subsection{Horton-Strahler ordering}
\label{hs}
\begin{figure}[p] 
\centering\includegraphics[width=.5\textwidth]{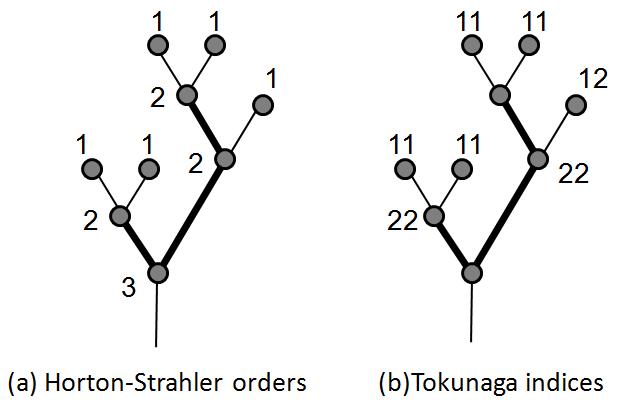}
\caption[Example of Horton-Strahler and Tokunaga indexing]
{Example of (a) Horton-Strahler ordering, and of (b) Tokunaga indexing.
Two order-2 branches are depicted by heavy lines in both panels.
The Horton-Strahler orders refer, interchangeably, 
to the tree nodes or to their parent links.
The Tokunaga indices refer to entire branches, and not to individual links.}
\label{fig_HST}
\end{figure}

The Horton-Strahler (HS) ordering of the vertices of a finite 
rooted labeled binary tree is performed in a hierarchical fashion, 
from leaves to the root \cite{NTG97,BWW00,Horton45,Strahler}: 
(i) each leaf has order $r({\rm leaf})=1$; 
(ii) when both children, $c_1, c_2$, of a parent vertex $p$ have 
the same order $r$, the vertex $p$ is assigned order $r(p)=r+1$; 
(iii) when two children of vertex $p$ have different orders, the vertex $p$ 
is assigned the higher order of the two.
Figure~\ref{fig_HST}(a) illustrates this definition.
Formally,
\begin{equation}
r(p)=\left\{
\begin{array}{ll}
r(c_1)+1&{\rm if~}r(c_1)=r(c_2),\\
\max\left(r(c_1),r(c_2)\right)&{\rm if~}r(c_1)\ne r(c_2).
\end{array}
\right.
\label{HS}
\end{equation}

A {\it branch} is defined as a union of connected vertices with 
the same order.
The branch vertex nearest to the root is called the 
{\it initial vertex}, the vertex farthest from the root 
is called the {\it terminal vertex}.
The order $\Omega(T)$ of a finite tree $T$ is the order 
$r(\phi)$ of its root, or, equivalently, the maximal order of 
its branches (or nodes).
The {\it magnitude} $m_i$ of a branch $i$ is the number of the 
leaves descendant from its initial vertex.
Let $N_r$ denote the total number of branches of order $r$
and $M_r$ the average magnitude of branches of order $r$ in
a finite tree $T$.

\begin{figure}[p] 
\centering\includegraphics[width=\textwidth]{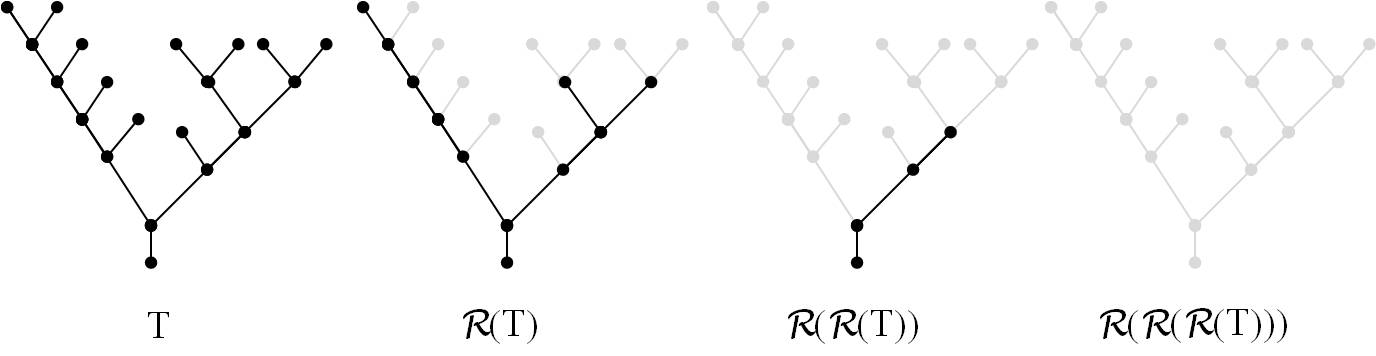}
\caption[Example of pruning]
{Example of consecutive application of the {\it pruning} operation
$\cR(\cdot)$ to the tree $T$.
In this example the tree has order $\Omega=3$ so $\cR^{(3)}(T)=\phi$.
For visual convenience the pruned branches are shown in
all panels by a light color.
Notice that pruning may produce chains of single-child nodes.}
\label{fig_pruning}
\end{figure}

An equivalent, and intuitively more appealing, definition of
the Horton-Strahler orders is done via the operation 
of {\it pruning} \cite{BWW00,Pec95}.
The pruning of an empty tree results in an empty tree,
$\cR(\phi)=\phi$.
The pruning $\cR(T)$ of a non-empty tree $T$, not necessarily binary, 
cuts the leaves and possible chains of degree-2 vertices 
connected to the leaves.
A vertex of degree 2 (or a single-child vertex) $v$ is defined 
by the conditions
$\langle v,1\rangle\in T$, 
$\langle v,2\rangle{\notin} T$.
Each chain of degree-2 vertices connected to a leaf
is uniquely identified by a vertex $v$ such that
$\langle v,u \rangle\in T$
implies 
$u=\langle 1,\dots, 1\rangle$.
The pruning operation is illustrated in Fig.~\ref{fig_pruning}.

The first application of pruning to a binary tree $T$
simply cuts the leaves, possibly producing some 
single-child vertices.
Some of those vertices are connected to the leaves via other
single-child vertices and thus will be cut at 
the next pruning, while the other occur deeper within the pruned 
tree and will wait for their turn to be removed.
It is readily seen that repetitive application of pruning
to any tree will result in the empty tree $\phi$.
The minimal $\Omega$ such that $\cR^{(\Omega)}(T)=\phi$
is called the {\it order} of the tree.  
A vertex $v$ of tree $T$ has the order $r$ if it has been 
removed at the $r$-th application of pruning:
$v\in\cR^{(k)}(T)\,\forall 1\le k <r$,
$v\notin\cR^{(r)}(T)$.
We say that a binary tree $T$ is {\it complete} if 
any of the following equivalent statements hold:
(i) each branch of $T$ consists of a single vertex;
(ii) orders of siblings (vertices with the common parent) 
are equal;
(iii) the parent vertex's rank is a unit higher than
that of each of its children.
There exists only one complete binary tree on
$n=2^k$ leaves for each $k=0,1,\dots$; all other trees are
called {\it incomplete}.

\subsection{Tokunaga indexing}
\label{tokunaga}
The Tokunaga indexing \cite{NTG97,Tok78,Pec95}
extends upon the Horton-Strahler orders; 
it is illustrated in Fig.~\ref{fig_HST}b.
This indexing focuses on incomplete trees by cataloging 
{\it side-branching}, which is the merging between branches 
of different order.
Let $\tau^k_{ij}$, $1\le k \le N_j$, $1\le i<j \le \Omega$ 
denotes the number of branches of order $i$
that join the non-terminal vertices of the $k$-th branch of 
order $j$. 
Then $N_{ij}=\sum_k\,\tau^k_{ij}$, $j>i$ is the total number 
of such branches in a tree $T$.
The Tokunaga index $T_{ij}$ is the average number of branches 
of order $i<j$ per branch of order $j$ in a finite tree of order 
$\Omega\ge j$:
\begin{equation}
T_{ij}= \frac{N_{ij}}{N_{j}}.
\label{tok}
\end{equation}

In a probabilistic set-up, one considers a space of finite binary 
trees with some probability measure.
Then, $N_i$, $\tau_{ij}^k$, $N_{ij}$, and $T_{ij}$ become random
variables.
We notice that if, for a given $\{ij\}$, the side-branch counts 
$\tau_{ij}^k$ are independent identically distributed random variables, 
$\tau_{ij}^k\stackrel{d}{=}\tau_{ij}$, then, by the law of large numbers,
\[T_{ij}\stackrel{\rm a.s}{\longrightarrow} {\sf E}\left(\tau_{ij}\right)
\quad {\rm as~} N_j\stackrel{\rm a.s}{\longrightarrow}\infty,\]
where the {\it almost sure} convergence 
$X_r \stackrel{\rm a.s.}{\longrightarrow} \mu$  
is understood as
${\sf P}\left(\displaystyle\lim_{r\to\infty}X_r=\mu\right)=1$.

For consistency, we denote the total number of order-$i$ branches
that merge with other order-$i$ branches by $N_{ii}$ and notice that
in a binary tree $N_{ii}=2\,N_{i+1}$.
This allows us to formally introduce the additional Tokunaga indices:
$T_{ii}=N_{ii}/N_{i+1}\equiv 2.$
The set $\{T_{ij}\}$, $1\le i \le \Omega-1,1\le j \le \Omega$, $i\le j$ 
of Tokunaga indices provides a complete statistical description of 
the branching structure of a finite tree of order $\Omega$.

\vspace{0.5cm}
Next, we define several types of tree self-similarity
based on the Horton-Strahler and Tokunaga indexing schemes.

\subsection{Horton laws}
\label{sst}

The {\it Horton laws}, widely observed in hydrological and 
biological networks \cite{TPN98,Horton45,VG00,DR00}, 
state, in their ultimate form,
\[
\frac{N_r}{N_{r+1}} = R_B,\quad
\frac{M_{r+1}}{M_r} = R_M,\quad
R_B,R_M>0,\quad r\ge 1,\]
where $N_r$, $M_r$ is, respectively, the total number and 
average mass of branches of order $r$ in a finite tree of 
order $\Omega$.
McConnell and Gupta \cite{MG08} emphasized the approximate,
asymptotic nature of the above empirical statements.
In the present set-up, it will be natural to formulate
the {\it Horton laws} as the almost sure convergence of 
the ratios of the branch statistics as the tree order 
increases:
\begin{eqnarray}
\label{NHL}
\frac{N_r}{N_{r+1}} &\stackrel{\rm a.s.}{\longrightarrow}& R_B>0,\quad
{\rm for~}r\ge 1, {~\rm as~}\Omega\to\infty,\\
\label{MHL}
\frac{M_{r+1}}{M_r} &\stackrel{\rm a.s.}{\longrightarrow}& R_M>0,\quad
{\rm as~}r,\Omega\to\infty.
\end{eqnarray}
Notice that the convergence in \eqref{NHL} is seen for 
the small-order branches, while the convergence in \eqref{MHL} 
--- for large-order branches.  
We call \eqref{NHL},\eqref{MHL} the {\it weak Horton laws}.
We also consider {\it strong Horton laws} that assume
an almost sure exponential dependence of the branch
characteristics on $r$ in a tree of finite order $\Omega$ and
magnitude $N$: 
\begin{eqnarray}
\label{NHLs}
N_r &\stackrel{\rm a.s.}{\sim}& N_0\,N\,R_B^{-r},\quad
{\rm for~}r\ge 1, {~\rm as~}\Omega\to\infty,\\
\label{MHLs}
M_r &\stackrel{\rm a.s.}{\sim}& M_0\,R_M^r,\quad
{\rm as~}r,\Omega\to\infty
\end{eqnarray}
for some positive constants $N_0,M_0,R_B$ and $R_M$ 
and with $x_r \stackrel{\rm a.s.}{\sim} y_r$ staying for 
\[{\sf P}\left(\displaystyle\lim_{r\to\infty} x_r/y_r =1\right)=1.\]
Clearly, the strong Horton laws imply the weak Horton laws.
The inverse in general is not true;
this can be illustrated by a sequence $M_r = R_M^r\,r^C$,
for any $C>0$, for which the weak Horton law \eqref{MHL} holds, 
while the strong law \eqref{MHLs} fails.
We notice also that $\Omega\to \infty$ implies $N\to\infty$,
but not vice versa; an example is given by a {\it comb} --- 
a tree of order $\Omega=2$ with an arbitrary number of side 
branches with Tokunaga index $\{12\}$.
This is why the limits above are taken with respect to $\Omega$,
not $N$. 
 
The strong Horton laws imply, in particular, that
\be
\label{ss}
N_r\stackrel{\rm a.s.}{\sim} const\,M_r^{-\alpha},
\quad\alpha=\frac{\log R_B}{\log R_M}
\ee
for appropriately chosen $r\to\infty$ and $\Omega\to\infty$,
for instance $r=\sqrt{\Omega}$.
The relationship \eqref{ss} is the simplest indication of 
self-similarity, as it connects the number $N_r$ and the size 
$M_r$ of branches via a power law.
However, a more restrictive property is conventionally required 
to call a tree self-similar; it is discussed in the next section.

\subsection{Tokunaga self-similarity}
In a deterministic setting, we call a tree $T$ of order $\Omega$  
a {\it self-similar tree} (SST) if its side-branching structure 
(i) is the same for all branches of a given order:
\[\tau_{ij}^k=:\tau_{ij},\quad 1\le k\le N_j,~1\le i<j\le\Omega,\]
and  
(ii) is invariant with respect to the branch order:
\be
\label{ss_d}
\tau_{i(i+k)}\equiv T_{i(i+k)} =: T_{k}\quad 
{\rm for~} 2\le i+k\le \Omega.
\ee
A {\it Tokunaga self-similar tree} (TSST) obeys an additional 
constraint first considered by Tokunaga \cite{Tok78}:
\be
\label{TSS}
T_{k+1}/T_k=c\quad \Leftrightarrow \quad
T_{k}=a\,c^{k-1}\quad a,c > 0,~1\le k\le\Omega-1.
\ee
In a random setting, we say that a tree $T$ of order $\Omega$ 
is self-similar if
${\sf E}\left(\tau_{i(i+k)}^j\right)=:T_k$
for $1\le j\le N_{i+k}$, $2\le i+k \le\Omega$; 
and it is Tokunaga self-similar if, 
furthermore, the condition \eqref{TSS} holds.

In a deterministic setting, for a tree satisfying
the weak Horton and Tokunaga laws\footnote{In a 
deterministic setting, the convergence in 
the Horton laws is understood as the convergence of 
sequences.}, one has \cite{Tok78,Pec95}:
\be
\label{HT}
R_B = \frac{2+c+a+\sqrt{(2+c+a)^2-8c}}{2}.
\ee
Peckham \cite{Pec95} has noticed that in a Tokunaga tree
of order $\Omega$ one has
$N_r=M_{\Omega-r+1}$, which implies that the Horton laws 
for masses $M_r$ follow from the Horton laws for the counts 
$N_r$ and $R_M=R_B$.
McConnell and Gupta \cite{MG08} have shown that the 
weak Horton laws with $R_B=R_M$ hold in a 
self-similar Tokunaga tree.
Zaliapin \cite{Zal10} has shown, moreover, that strong 
Horton laws hold in a Tokunaga tree and, at the same time, 
even weak Horton laws may not hold in a general, non-Tokunaga, 
self-similar tree. 

The Tokunaga self-similarity describes a two-parametric class
of trees, specified by the Tokunaga parameters $(a,c)$.
Our goal is to demonstrate that the Tokunaga class is
not only structurally simple but is also sufficiently wide.
This study establishes the Tokunaga self-similarity for
the level-set trees of symmetric homogeneous Markov chains, 
and, as a direct consequence, for the trees of their scaling limits 
including a regular Brownian motion.

\subsection{Stochastic self-similarity}
Burd {\it et al.} \cite{BWW00} define {\it stochastic self-similarity}
for a random tree $\tau\in\left(\mathbb{T}_0,P\right)$
as the distributional invariance with respect to the pruning $\cR(\tau)$:
\[P\left(\cdot|\tau\ne\phi\right)\circ \cR^{-1}=P(\cdot )\]
and prove the following result that explains the importance of 
Tokunaga self-similarity within the class of Galton-Watson trees 
as well as the special role of the Galton-Watson critical binary trees.
\begin{thm}{\rm \cite[Theorems 1.1, 1.2, 3.17]{BWW00}}
Let $\tau\in\left(\mathbb{T}_0,GW_{\{p_k\}}\right)$
with bounded offspring number.
Then the following statements are equivalent:
\begin{itemize}
\item[(i)] Tree $\tau$ is stochastically self-similar.
\item[(ii)] ${\sf E}(\tau_{i(i+k)}) =: T_k$, {\it i.e.},
the expectation is a function of $k$ and $T_k$ is defined by this equation.
\item[(iii)] Tree $\tau$ has the critical binary offspring distribution,
$p_0=p_2=1/2$.
\end{itemize}
\end{thm} 

These authors show, furthermore, how the arbitrary binary Galton-Watson 
distribution is transformed under the operation of pruning.
\begin{thm}{\rm \cite[Proposition 2.1]{BWW00}}
\label{tree_prune}
Let $\tau$ be a finite tree with a binary Galton-Watson distribution,
$p_0+p_2=1$, with $p_2\le 1/2$. 
Let $\tau_{n+1}=\cR(\tau_n)$, $n\ge 0$, $\tau_0=\tau$.
Then $\tau_{n+1}$ has the binary Galton-Watson distribution
$p_0^{(n+1)}+p_2^{(n+1)}=1$ with
\[p_2^{(n+1)}=\frac{\left[p_2^{(n)}\right]^2}
{\left[p_0^{(n)}\right]^2+\left[p_2^{(n)}\right]^2}.\]
\end{thm}

We demonstrate below that stochastic (or distributional) 
self-similarity, within the class of tree representations of 
homogeneous Markov chains, holds only for Markov
chains with symmetric exponential increments.

\section{Main results}
\label{main}
Let $X_k$, $k\in\mathbb{Z}$ be a real valued Markov chain with homogeneous
transition kernel $K(x,y)\equiv K(x-y)$, for any $x,y\in\mathbb{R}$.
We call $X_k$ a homogeneous Markov chain (HMC).
When working with trees, $X_k$ will also denote a function from 
$C(\mathbb{R})$ obtained by liner interpolation of the values of 
the original time series $X_k$; this create no ambiguities in
the present context.
 
A HMC is called {\it symmetric} (SHMC) if its transition kernel
satisfies $K(x)=K(-x)$ for any $x\in\mathbb{R}$.
We call an HMC {\it exponential} (EHMC) if its kernel
is a mixture of exponential jumps.
Namely,
\[K(x)=p\,\phi_{\lambda_u}(x)+(1-p)\,\phi_{\lambda_d}(-x),
\quad 0\le p \le 1, \lambda_u,\lambda_d>0,\]
where $\phi_{\lambda}$ is the exponential density
\be
\label{exp}
\phi_{\lambda}(x)=
\begin{cases} 
\lambda e^{-\lambda x}, & x \geq 0, \\ 
0, & x<0. 
\end{cases}
\ee
We will refer to an EHMC by its parameter triplet $\{p,\lambda_u,\lambda_d\}$.

The concept of tree self-similarity is based on the notion of
{\it branch order} and is tightly connected to the {\it pruning}
operation (Sect.~\ref{hs}, Fig.~\ref{fig_pruning}).
In terms of time series (or tamed real functions), pruning corresponds 
to coarsening the time series resolution by removing the 
local maxima.
An iterative pruning corresponds to iterative transition to
the local minima. 
We formulate this observation in the following proposition.

\begin{prop}
\label{ts_prune}
The transition from a time series $X_k$ to the time series 
$X^{(1)}_k$ of its local minima corresponds to the pruning 
of the level-set tree $\textsc{level}(X)$.
Formally,
\[\textsc{level}\left(X^{(m)}\right) 
= \cR^m\left(\textsc{level}(X)\right), \forall m \ge 1,\]
where $X^{(m)}$ is obtained from $X$ by iteratively taking
local minima $m$ times (i.e., local minima of local minima
and so on.)
\end{prop}

The next result establishes invariance of several classes
of Markov chains with respect to the pruning operation.

\begin{lem}
\label{inv_prune}
(a) The local minima of a HMC form a HMC.
(b) The local minima of a SHMC form a SHMC.
(c) The local minima of an EHMC with parameters 
$\{p,\lambda_u,\lambda_d\}$ form a EHMC with parameters 
$\{p^*,\lambda_u^*,\lambda_d^*\}$, where
\be 
\label{iteration}
p^*=\frac{p\,\lambda_d}{p\,\lambda_d+(1-p)\,\lambda_u}, \quad 
\lambda^*_d=p\lambda_d,~~\text{ and }~~
\lambda^*_u=(1-p)\lambda_u.
\ee
\end{lem}

Let $\{M_t\}\equiv \{M^{(1)}_t\}$, $t\in \mathcal{T}_1\subset\mathbb{R}$, 
be the set of local minima of $X_t$, not including the boundary minima;
$\{M^{(2)}_t\}$, $t\in \mathcal{T}_2\subset\mathbb{R}$, 
be the set of local minima of local minima
(local minima of second order), {\it etc.}, with $\{M^{(j)}_t\}$, 
$t\in \mathcal{T}_j\subset\mathbb{R}$ being the local minima of order $j$.
We call a segment between two consecutive points from $\mathcal{T}_r$, 
$r\ge1$, a {\it (complete) basin} of order $r$.
For each $r$, there might exist a single leftmost and 
a single rightmost segments of $X_t$ that do not belong 
to any basin or order $r$, with a possibility
for them to merge if $X_t$ does not have basins
of order $r$ at all. 
We call those segments {\it incomplete} basins of order $r$.
There is a bijection between basins (complete and incomplete) 
of order $r$ in $X_t$ and branches of Horton-Strahler 
order $r$ in $\textsc{level}(X_t)$.
This explains the terms {\it complete branch} and 
{\it incomplete branch} of order $r$.

\begin{thm}[{\bf Horton and Tokunaga self-similarity}]
\label{T1}
The combinatorial level set tree 
$\textsc{shape}\left(\textsc{level}(X),1\right)$ 
of a finite SHMC $X_k$, $k=1,\dots,N$ 
satisfies the strong Horton laws for any $r\ge 1$, 
asymptotically in $N$:
\be
\label{Horton}
N_r\overset{\rm a.s.}{\sim} N\,R_{\rm B}^{-r},\quad R_{\rm B}=4,
\quad{\rm as~}N\to\infty.
\ee 
Furthermore, $T=\textsc{shape}\left(\textsc{level}(X,1)\right)$
is a Tokunaga self-similar tree with parameters $(a,c)=(1,2)$.
Specifically, for a finite tree $T$ of order $\Omega(N)$
the side-branch counts $\tau_{i(i+k)}^j$ with
$2\le i+k\le \Omega$
for different complete branches $j$ of order $(i+k)$ 
are independent identically distributed random variables such that 
$\tau_{i(i+k)}^j\overset{d}{=:}\tau_{i(i+k)}$
and
\be
\label{Tokunaga}
{\sf E}\left[\tau_{i(i+k)}\right]=:T_k = 2^{k-1}.
\ee
Moreover, $\Omega\overset{a.s.}{\to}\infty$ as
$N\to\infty$ and, for any $i,k\ge 1$, we have
\[T_{i(i+k)}\overset{\rm a.s.}{\longrightarrow} T_k=2^{k-1},
\quad{\rm as~}N\to\infty,\]
where $T_{i(i+k)}$ can be computed over the entire $X_k$.
\end{thm}

\begin{figure}[p] 
\centering\includegraphics[width=0.6\textwidth]{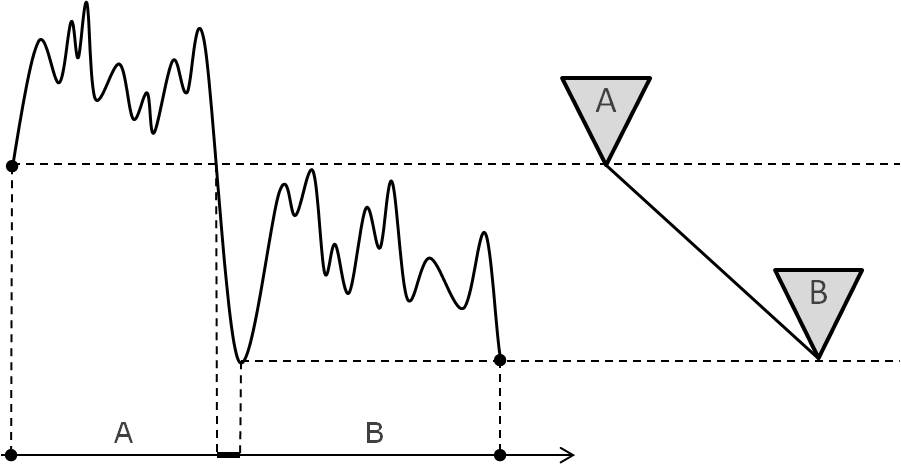}
\caption{Illustration of tree construction for an infinite time series.
The time series $X_t$ is divided here into two vertically shifted excursions,
marked {\bf A} and {\bf B} in the time axis, and one fall,
depicted by the heavy segment on the time axis. 
The descending ladder $L_X$ consists of two isolated 
points and one interval (heavy segment on the time axis).
The excursions correspond to the two trees represented by marked 
triangles, the interval from the descending ladder corresponds 
to the line that connects the trees A and B. }
\label{fig_ladder}
\end{figure}

Next we extend this result to the case of infinite time 
series and the weak limits of finite time series.
For a linearly interpolated time series  $X_t$, $t\ge 0$ 
(equivalently, for a continuous function with a countable 
number of separated local extrema) consider the 
{\it descending ladder} $L_X=\{t:X_t=\underline{X}[0,t]\}$,
which in our settings is a set of isolated points and 
non-overlapping intervals (Fig.~\ref{fig_ladder}).
The function $X_t$ is naturally divided into a series
of vertically shifted positive excursions on the
intervals not included in $L_X$ and monotone falls on 
the intervals from $L_X$. 
Any (in the a.s. sense) infinite SHMC can be decomposed
into infinite number of such finite excursions and 
finite falls.   
We will index the excursions by index $i\ge 1$ from left to right.
The extreme time series $\cE\left(X_k^i\right)$ for 
each finite excursion $X_t^{i}$ is a Harris path
for a finite tree $\textsc{level}\left(X_t^{i}\right)$.
Hence, each such finite excursion completely specifies a single
subtree of $\textsc{tree}\left(X_t\right)$.
In particular, it completely specifies the HS orders 
for all vertices and 
Tokunaga indices for all branches except the one 
containing the root within 
$\textsc{level}\left(X_t^{i}\right)$. 
We also notice that each fall of $X_t$ on an interval 
from $L_X$ corresponds
to an individual edge of $\textsc{tree}\left(X_t\right)$.
Combining the above observations, we conclude that 
the tree
$\textsc{tree}\left(X_t\right)$
can be represented as infinite number of 
subtrees $\textsc{level}(X_t^i)$ connected
by edges that correspond to the falls
of $X_t$ on the descending ladder, see Fig.~\ref{fig_ladder}.
Pitman calls this construction, applied to the standard
Brownian motion rather than time series, {\it a forest of
trees attached to the floor line} \cite[Section 7.4]{Pitman}.
Let $N_r^n$ and $N_{ij}^n$ denote, respectively, the number 
of branches of order $r$ and the number of side branches 
of Tokunaga index $\{ij\}$ in the first $n$ excursions
of $X_t$ as described above.
We introduce the cumulative quantities
\[\eta_r^n:=\frac{N_r^n}{N_{r+1}^n},\quad 
T_{ij}^n:=\frac{N_{ij}^n}{N_j^n}\]
and define, for the infinite time series $X_t$, 
\be
\label{inf_ss}
\eta_r(X_t) = \lim_{n\to\infty} \eta_r^n,\quad
T_{ij}(X_t) = \lim_{n\to\infty} T_{ij}^n,
\ee
whenever the above limits exist in an appropriate 
probabilistic sense.

By Proposition~\ref{inv1}, the level set tree of a finite 
excursion $X_t^k$ is not affected by monotonic transformations
of time and value. 
This allows to expand the above definition \eqref{inf_ss} 
to the weak limits of time series via the the Donsker's theorem.
In particular, if $X_t$ is a SHMC whose increments have standard
deviation $\sigma$, then the rescaled segments $X_t$ weakly 
converge to the regular Brownian motion $B_t$, $0\le t\le 1$.
Namely,
\[X_{(nt)}/\sqrt{n}\overset{d}{\to} \sigma\,B_t\] 
as $n\to\infty$
through the end point of the finite excursions that
comprise $X_t$.
This leads to the following result. 
\begin{cor}
\label{Brown}
The combinatorial tree 
$\textsc{shape}\left(\textsc{tree}(B_t)\right)$ of a 
regular Brownian motion $B_t$, $t\in[0,1]$
satisfies the Horton and Tokunaga self-similarity laws.
Namely,
\be
\eta_r(B_t) = 4 ~{\rm for~} r\ge 1\quad{\rm and}\quad
T_{i(i+k)}(B_t)=2^{k-1}~{\rm for~} i,k\ge 1,
\ee
where the limits \eqref{inf_ss} are understood in the almost sure sense.
\end{cor}

We conclude this section with a conjecture motivated by the above
result as well as extensive numeric simulations \cite{Webb09}.
\begin{conj}
\label{BH}
The tree $\textsc{shape}\left(\textsc{tree}\left(B^H\right)\right)$ of 
a fractional Brownian motion $B^H_t$, 
$t\in[0,1]$ with the Hurst index $0<H<1$ is Tokunaga self-similar with
$T_{i(i+k)}(B^H)=T_k=c^{k-1}$, $c=2H+1$, $i,k\ge 1$. 
According to \eqref{HT}, this corresponds to the Horton self-similarity with
\be
\eta_r(B^H) = 2+H+\sqrt{H^2+2}, r\ge 1.
\ee
The sense of limits \eqref{inf_ss} is to be determined.
\end{conj}

\section{Exponential chains}
\label{add}
This section focuses on exponential chains, which enjoy
an important distributional self-similarity and whose 
level-set trees have the Galton-Watson distribution. 

\subsection{Distributional self-similarity} 
\label{DSS_sec}
Consider a SHMC $X_k$, $k\in\mathbb{Z}$ with kernel
\[K(x)=\frac{f(x)+f(-x)}{2},\]
where $f(x)$ is a probability density function with support $\mathbb{R}^+$. 
The series of local minima of $X_k$ (or, equivalently, pruning $X_k^{(1)}$ of $X_k$) 
also forms a SHMC with transition kernel $K_1(x)$ 
(see Lemma~\ref{inv_prune}(b)).
It is natural to look for chains invariant with respect to
the pruning:
\be
\label{d_inv}
X_k \overset{d}{=} c\,X_k^{(1)}, c>0.
\ee
By Proposition~\ref{inv1}, such invariance would guarantee
the {\it distributional} Tokunaga self-similarity:
\be
\label{Tok_d}
\tau_{i(i+k)}^j \overset{d}{=:} \tau_{i(i+k)}=T_k,
\quad 1\le j\le N_{i+k}, 1\le i+k\le\Omega,
\ee
where $T_k$ is a random number of side-branches of order $i$
that join an arbitrarily chosen branch of order $(i+k)$.
Hence, we seek the conditions on $f(x)$ to ensure that 
$K_1(x)=c^{-1}K(x/c)$
for some constant $c>0$. 

\begin{prop}
\label{DSS}
The local minima of a SHMC $X_k$ with kernel $K(x)$ form a SHMC with kernel
\[K_1(x)=\frac{K(x/c)}{c},\quad c>0\]
if and only if $c=2$ and
\be
\label{laplace}
\Re\left[\widehat{f}(2s)\right]=\left|\frac{\widehat{f}(s)}{2-\widehat{f}(s)}\right|^2,
\ee
where $\widehat{f}(s)$ is the characteristic function of $f(x)$ and
$\Re[z]$ stays for the real part of $z\in\mathbb{C}$.
\end{prop}

Observe that the set of densities $f(x)$ that satisfy (\ref{laplace}) is not empty. 
A solution is given for example by the Laplace density with $\lambda>0$,
{\it e.g.} for $f(x)=\phi_{\lambda}(x)$ with exponential density $\phi_{\lambda}(x)$
of \eqref{exp}, that is by an EHMC $\{1/2,\lambda,\lambda\}$.

\subsection{Distributional self-similarity for symmetric exponential chains}
\label{expj}
Lemma~\ref{inv_prune}(c) allows one to study the behavior of the 
EHMCs formed by local minima, minima of minima, and so on of
an EHMC $X_k$ with parameters $\{p,\lambda_u,\lambda_d\}$.
Introducing the variables
\be
\label{ag}
A=\frac{1-p}{p},\quad \gamma = \frac{\lambda_d}{\lambda_u}
\ee
one readily obtains that their counterparts $\{A^*,\gamma^*\}$ for 
the chain of local minima, given by \eqref{iteration}, are expressed as
\be
\label{dyn}
A^* = \frac{A}{\gamma},\quad \gamma^*=\frac{\gamma}{A}.
\ee
Notably, this means that the chain of local minima for 
{\it any} EHMC form an EHMC with $A\,\gamma=1$.
The only fixed point in the space $(A,\gamma)$ with 
iteration rules \eqref{dyn} is the point $(A=1,\gamma=1)$, 
which corresponds to the distributionally self-similar
EHMS discussed in Sect.~\ref{DSS_sec}.
This point is an image (under the pruning operation) of the
EHMCs with $A=\gamma$ or 
$p\,\lambda_d = (1-p)\,\lambda_u$.
The last condition is equivalent to ${\sf E}(X_k-X_{k-1}) = 0$
for any $k>1$.
The chain of local minima for any EHMC with $A>\gamma$ ($A<\gamma$)  
corresponds to a point on the upper (lower) part of the hyperbola 
$A\,\gamma=1$.
Any point on this hyperbola, except the fixed point $(1,1)$, moves
away from the fixed point toward $(0,\infty)$ or $(\infty,0)$.
This is illustrated in Fig.~\ref{fig4}.
It follows that the Tokunaga and even weaker Horton self-similarity 
is only seen for a symmetric EHMC.
The above discussion can be summarized in the following statement. 
\begin{thm}
\label{eH}
Let $X_k$ be an EHMC $\{p,\lambda_u,\lambda_d\}$.  
Then $X_k$ satisfies the distributional self-similarity \eqref{d_inv}
if and only if $p=1/2$, $\lambda_u=\lambda_d$.
Furthermore, the multiple pruning $X_k^{(m)}$, \mbox{$m>1$} of $X_k$ 
satisfies the distributional self-similarity \eqref{d_inv}
if and only if the chain's increments have zero mean, or, equivalently,
if and only if $p\,\lambda_d = (1-p)\,\lambda_u$.
In this case, the self-similarity is achieved after
the first pruning, that is for the chain $X_k^{(1)}$ of 
local minima. 
\end{thm}

\begin{cor}
The regular Brownian motion with drift is not Tokunaga self-similar. 
\end{cor}

\subsection{Connection to Galton-Watson trees}
An important, and well known, fact is that the Galton-Watson 
distribution (see Sect.~\ref{GW}) is the 
characteristic property of trees that have Harris paths 
with alternating exponential steps.
We formulate this result using the terminology of our paper.
\begin{thm}{\rm \cite[Lemma 7.3]{Pitman},\cite{LeGall93,NP89}}
\label{Pit7_3}
Let $X_k$ be a discrete-time excursion with finite 
number of local minima.
The level set tree $\textsc{shape}\left(\textsc{level}(X_k,1)\right)$ 
is a binary Galton-Watson tree with $p_0+p_2=1$ if and only if the 
rises and falls of $X_k$, excluding the last fall, are 
distributed as independent exponential 
variables with parameters $(\mu+\lambda)$
and $(\mu-\lambda)$, respectively, for $0\le \lambda < \mu$.
In this case, 
\[p_0=\frac{\mu+\lambda}{2\,\mu},\quad
p_2=\frac{\mu-\lambda}{2\,\mu}.\]
\end{thm}

We now use this result to relate sequential pruning 
of Galton-Watson trees (see Theorem~\ref{tree_prune}) and
pruning of EHMCs.
Consider the first positive excursion $X_k$ of an 
EHMC with parameters 
$\{p^{(0)}= p=1-q,\lambda_u,\lambda_d\}$.
The geometric stability of the exponential distribution 
implies that the monotone rises and falls of $X_k$
are exponentially distributed with parameters
$q\,\lambda_u$ and $p\,\lambda_d$, respectively.
The Theorem~\ref{Pit7_3} implies that 
$\textsc{shape}\left(\textsc{level}(X_k)\right)$ is
distributed as a binary Galton-Watson tree,
$p_0+p_2=1$, with 
\be
\label{p20}
p_2 \equiv p_2^{(0)}=\frac{p\,\lambda_d}{q\,\lambda_u+p\,\lambda_d}.
\ee
The first pruning $X_k^{(1)}$ of $X_k$, according to \eqref{iteration},
is the EHMC with parameters
\[\left\{p^{(1)}=\frac{p\,\lambda_d}{q\,\lambda_u+p\,\lambda_d},
q\,\lambda_u,p\,\lambda_d\right\}.\]
Its upward and downward monotone increments are exponentially 
distributed with parameters, respectively,
\[\frac{(q\,\lambda_u)^2}{q\,\lambda_u+p\,\lambda_d}\quad
{\rm and}\quad
\frac{(p\,\lambda_d)^2}{q\,\lambda_u+p\,\lambda_d}.\]
By Theorem~\ref{Pit7_3}, the level-set tree for an arbitrary positive 
excursion of $X_k^{(1)}$ is a binary Galton-Watson tree,
$p_0^{(1)}+p_2^{(1)}=1$, with 
\[p_2^{(1)}=\frac{(p\,\lambda_d)^2}{(q\,\lambda_u)^2+(p\,\lambda_d)^2}.\]

Continuing this way, we find that $n$-th pruning $X_k^{(n)}$ of 
$X_k\equiv X_k^{(0)}$ is an EHMCs such that the level set tree of 
its arbitrary positive excursion have a binary Galton-Watson 
distribution, $p_0^{(n)}+p_2^{(n)}=1$, with
\[p_2^{(n)}=\frac{(p\,\lambda_d)^{2^n}}
{(q\,\lambda_u)^{2^n}+(p\,\lambda_d)^{2^n}}.\]
This can be rewritten in recursive form as
\[p_2^{(n)}=\frac{\left[p_2^{(n-1)}\right]^2}
{\left[p_0^{(n-1)}\right]^2+\left[p_2^{(n-1)}\right]^2},\quad n\ge 1\]
with $p_2^{(0)}$ given by \eqref{p20}.
Notably, this is the same recursive system as that discovered
by Burd {\it et al.} \cite[Proposition 2.1]{BWW00} (see 
Theorem~\ref{tree_prune} above) in their analysis of 
consecutive pruning for the Galton-Watson trees.
Another noteworthy relation is given by
\[p^{(n)} = p_2^{(n-1)},\quad n\ge 1, \quad p^{(0)}=p, p_2^{(0)}=p_2,\]
which connects the ``horizontal'' probability $p^{(n)}$ of 
an upward jump in a pruned time series $X_k^{(n)}$ with the
``vertical'' probability $p_2^{(n-1)}$ of branching in
a Galton-Watson tree.


\section{Terminology and proofs}
\label{proofs}
\subsection{Level-set trees: Definitions and terminology}
\label{terms}
This section introduces terminology for discussing the 
hierarchical structure of the local extrema of a finite 
time series $X_k$ and relating it to the level 
set tree $\textsc{level}(X)$.
For consistency we repeat some terms introduced 
above to formulate Theorem~\ref{T1}.

Let $\{M_t\}\equiv \{M^{(1)}_t\}$, $t\in \mathcal{T}_1\subset\mathbb{R}$, 
be the set of local minima of $X_t$, not including possible boundary minima;
$\{M^{(2)}_t\}$, $t\in \mathcal{T}_2\subset\mathbb{R}$, 
be the set of local minima of local minima (local minima of second order), 
{\it etc.}, with $\{M^{(j)}_t\}$, $t\in \mathcal{T}_j\subset\mathbb{R}$ 
being the local minima of order $j$.
Next, let $\{m_s\}\equiv \{m^{(1)}_s\}$, $s\in\mathcal{S}_1\subset\mathbb{R}$, 
be the set of local maxima of $X_k$, including possible boundary maxima, and 
$\{m^{(j+1)}_s\}$, $s\in\mathcal{S}_{j+1}\subset\mathbb{R}$ the set of local 
maxima of $\{M^{(j)}_t\}$ for all $j\ge 1$.
We will call a segment between two consecutive points from $\mathcal{T}_j$ 
a {\it (complete) basin} of order $j$.
Clearly, $\mathcal{T}_1\supset \mathcal{T}_2 \supset\dots$ and each
basin of order $r$ is comprised of a non-zero number of basins
of arbitrary order $k<r$.  
For each $r$, there might exist a single leftmost and 
a single rightmost segments of $X_t$ that do not belong 
to any basin or order $r$, with a possibility
for them to merge if $X_t$ does not have basins
of order $r$ at all. 
We call those segments {\it incomplete} basins of order $r$.

By construction, each basin of order $j$ contains exactly one point from 
$\mathcal{S}_j$; e.g., there is a single local maximum 
from $\mathcal{S}_1$ between two consecutive local minima 
from $\mathcal{T}_1$, {\it etc.}
There exists a bijection between basins (complete and incomplete) 
of order $r$ in $X_t$ and branches of Horton-Strahler 
order $r$ in $\textsc{level}(X_t)$;
this explains the terms {\it complete branch} and 
{\it incomplete branch} of order $r$.
More specifically, there is a bijection between the terminal 
vertices of order-$r$ branches --- {\it i.e.,} vertices parental 
to two branches of order $(r-1)$ --- and the local maxima
from $\mathcal{S}_j$ within the respective basins.

Let us fix an arbitrary local minimum $X_k$ of order $r_k$; then
$k\in \mathcal{T}_j$ for $1\le j \le r_k$ and
$k\notin \mathcal{T}_j$ for $j>r_k$.
For each $j > r_k$ there exists a unique basin of order $j$ that contains $k$;
we denote the boundaries of this basin by 
$l_k^{(j)},r_k^{(j)}\in \mathcal{T}^{(j)}$, $l_k^{(j)}<r_k^{(j)}$.
Denote by $c^{(j)}_k$ the unique point from $\mathcal{S}_j$ within the 
interval $\left(l_k^{(j)},r_k^{(j)}\right)$.
Multiple points $X_k$ may correspond to the same triplet
$\left(l_k^{(j)},c^{(j)}_k,r_k^{(j)}\right)$, which will create no confusion.
These definitions are illustrated in Fig.~\ref{fig1}.

Consider now a point $k$ of local minimum such that
$k\notin\cup_{j\ge 1} m^{(j)}_s$.
If $l_k^{(j)} < k < c^{(j)}_k$ for a given $j>r_k$ then we call 
the point $l_k^{(j)}$ the local minimum of order $j$ {\it adjacent} to $k$ and
the point $r_k^{(j)}$ the local minimum of order $j$ {\it opposite} to $k$.
The analogous terminology is introduced in case $c_k^{(j)} < k < r^{(j)}_k$.
By construction, $X_k$ is always greater than the value of its
adjacent minimum of any order $j>r_k$.
The value of the opposite minimum of order $j$ is denoted by $M^{(j)}_k$.
We have, for each $k$,
\be
\label{morder}
M^{(1)}_k \ge M^{(2)}_k \ge M^{(3)}_k \ge \dots
\ee

We already noticed that the local maxima $m^{(1)}_t$ correspond to the 
tree leaves, that is to its branches of Horton order $r=1$.
The set $m^{(j+1)}_t$ for each $j\ge 1$  corresponds to 
the vertices parental to two branches of the same HS order $j$; 
they are the terminal vertices of order-$(j+1)$ branches.
All other local minima of $X_k$ correspond to vertices parental 
to two vertices of different SH order; we will refer to this as 
{\it side-branching}.
Specifically, a local minimum $X_k$ of order $i$ forms a side-branch of 
order $\{ij\}$ if
\be
\label{sb}
M^{(j-1)}_k \ge X_k \ge M^{(j)}_k,
\ee
where the first inequality disappears when $j=i+1$.
Figure~\ref{fig2} illustrates this for a basin of second order. 
In general, each basin of order $r$ contains a uniquely 
specified positive excursion attached to its higher end.
The local maxima of order $k<r$ from this excursion
correspond to the side-branches with Tokunaga index $\{km\}$
with $m\le r$.
The local maxima of order $k<r$ within the basin but outside 
of this excursion     
correspond to the side-branches with Tokunaga index $\{km\}$
with $m> r$.

\subsection{Proofs}

\noindent
{\bf Proof of Propositions~\ref{inv1},\ref{inv2},\ref{Harris} and \ref{ts_prune}:}
The statements readily follow from the definition of level set trees.
\qed

\vspace{.5cm}
\noindent
{\bf Proof of Lemma~\ref{inv_prune}:}

(a) Follows from the independence of increments in $X_k$.

(b)
Let $\{M_j\}$ be the sequence of local minima of $X_k$ and
$d_j = M_{j+1}-M_j$.
We have, for each $j$
\be
\label{twosum}
d_j = \sum_{i=1}^{\xi_+} Y_i - \sum_{i=1}^{\xi_-} Z_i,
\ee
where $\xi_+$ and $\xi_-$ are independent geometric random variables 
with parameter $1/2$: 
\[{\sf P}(\xi_+=k) = {\sf P}(\xi_-=k) = 2^{-k}, \quad k=1,2,\dots;\]
$Y_i$, $Z_i$ are independent identically distributed (i.i.d.) random variables with density 
$f(x)$.
Here the first sum corresponds to $\xi_+$ positive increments of $X_k$
between a local minimum $M_j$ and the subsequent local maximum $m_j$ and
the second sum to $\xi_-$ negative increments between the local maximum
$m_j$ and the subsequent local minimum $M_{j+1}$.
It is readily seen that both the sums in \eqref{twosum} have the same 
distribution, and hence their difference has a symmetric distribution.
We notice that the symmetric kernel for the sequence of local minima $\{M_j\}$ 
is necessarily different from $K(x)$.

(c)
Consider an EHMC $X_k$ with parameters $\{p,\lambda_u,\lambda_d\}$.
By statement (a) of this lemma, 
the local minima of $X_k$ form a HMC with transition kernel $K_1(x)$.
The latter is the probability distribution of the jumps $d_j$
given by \eqref{twosum} with $\xi_+$, $\xi_-$ being geometric random variables
with parameters $p$ and $(1-p)$ respectively, 
$Y_i\overset{d}{=}\phi_{\lambda_u}$, and
$Z_i\overset{d}{=}\phi_{\lambda_d}$.
For the characteristic function of $K_1$ one readily has
\[\widehat{K_1}(s)=
\frac{p(1-p)\lambda_d \lambda_u}{\left((1-p)\lambda_u-is\right)(p\lambda_d+is)}
=p_* \cdot \widehat{\phi_{\lambda^*_u}}(s)+(1-p_*) 
\cdot \widehat{\phi_{\lambda^*_d}}(-s)\]
with
\[
p_*=\frac{p\,\lambda_d}{p\,\lambda_d+(1-p)\,\lambda_u}, \quad 
\lambda^*_d=p\lambda_d,~~\text{ and }~~
\lambda^*_u=(1-p)\lambda_u.
\]
Thus
\[K_1(x)=p_* \phi_{\lambda^*_u}(x)+(1-p_*)\phi_{\lambda^*_d}(-x).\]
This means that the HMC of local minima also jumps according to a two-sided 
exponential law, only with different parameters $p_*$, $\lambda^*_d$ and $\lambda^*_u$. 
\qed

\vspace{.5cm}
\noindent
{\bf Proof of Theorem~\ref{T1}: Horton self-similarity}

We notice that the number $N_r$ of order-$r$ branches 
in $\textsc{level}(X)$ equals the number 
$|\cS_r|$ of local maxima $m^{(r)}_s$ of order 
$r$ (with the convention that the local maxima 
of order 0 are the values of $X_k$).
The probability for a given point of $X_k$ to be a local
maximum equals the probability that this point is higher
than both its neighbors.
The Markov property and symmetry of the chain imply that
this probability is $1/4$.
Hence the average number of local maxima is
\[{\sf E}\left(|\cS_0|\right)
={\sf E}\left(N_1\right) 
= \sum_{i=2}^{N-1}{\sf P}(X_{i-1}<X_i>X_{i+1})
= \frac{N-2}{4}\sim \frac{N}{4}.\]
Let $l_i$ denote the event ($X_i$ is a local maximum).
By Markov property, the events $l_i$, $l_j$ are independent for
$|i-j|\ge 2$;
hence, the variance ${\sf V}(N_1)\propto N$.
This yields
\[\lim_{N\to\infty}{\sf E}\left(\frac{N_1}{N}\right)=1/4,\quad
\lim_{N\to\infty}{\sf V}\left(\frac{N_1}{N}\right)=0.\]
One can combine the strong laws of large numbers for (i) the 
proportion of the upward increments of $X_t$ (that converges
to 1/2) and (ii) the proportion of upward increments followed
by a downward increment (that converges to 1/2)
to obtain
$N_1/N \overset{\rm a.s.}{\to} 1/4$, and, in particular,
$N_1\overset{\rm a.s.}{\to}\infty$ as $N\to\infty$.

We use now Lemma~\ref{inv_prune}(b) to find, applying the same argument
to the pruned time series, that
$N_r/N_{r-1}\overset{\rm a.s.}{\to} 1/4$ as $N\to\infty$ for any $r> 1$.
Finally,
\[\frac{N_r}{N}=
\frac{N_r}{N_{r-1}}\frac{N_{r-1}}{N_{r-2}}\dots
\frac{N_{1}}{N} \overset{\rm a.s.}{\longrightarrow} 4^{-r},\quad N\to\infty,\]
which completes the proof of the strong Horton law \eqref{Horton}.
\qed

The proof of the Tokunaga self-similarity will require several 
auxiliary statements formulated below.

\begin{lem}
\label{L4}
A basin of order $j$ contains on average $4^{j-k}$ basins of 
order $k$, for any $j > k\ge 1$.
\end{lem}
\noindent
{\bf Proof of Lemma~\ref{L4}:}
We show first that a basin of order $(j+1)$ contains on average 
4 local minima of order $j\ge 1$.
The number $\xi$ of points of $X_k$ within a first-order basin
({\it i.e.}, between two consecutive local minima)
is $\xi=1 + \xi_++\xi_-$, where $\xi_+$, $\xi_-$ are, respectively, 
the numbers of basin points (excluding the basin boundaries) to the 
left and right of its local maximum $m$; and the latter is counted 
separately in the expression above.
The independence of increments of $X_k$ impies
\[{\sf P}(\xi_+=k)={\sf P}(\xi_-=k)=2^{-k-1}, k=0,1,\dots,\]
and hence
\be
\label{Exi}
{\sf E}[\xi]=1+{\sf E}[\xi_+]+{\sf E}[\xi_-]=1+1+1=3.
\ee
By Lemma~\ref{inv_prune}(b), the same result holds for the average number
of local minima of order $j$ within an order-$(j+1)$ basin,
for any $j\ge 1$.
Thus, the average number of order-$j$ basins within an order-$(j+1)$ basin
is  ${\sf E}[\xi]+1=4.$

The independence of increments of $X_k$ implies that the number
of order-$(j-1)$ subbasins within an order-$j$ basin is independent
of the numbers of order-$j$ basins within an order-$(j+1)$ basin.
This leads to the Lemma's statement.
\qed

\begin{lem}
\label{L2}
Let $a$ and $b$ be two points chosen at random and without replacement 
from the set $\{1,2,\dots,N\}$ and $\eta=(\eta_1,\eta_2,\eta_3)$ denotes 
the random number of points within the following intervals respectively:
(i)   $\left[1,\,\min(a,b)\right)$,
(ii)  $\left(\min(a,b),\,\max(a,b)\right)$, and
(iii) $\left(\max(a,b),\, N\right]$.
Then the triplet $\eta$ has an exchangeable distribution.
\end{lem}

\noindent{\bf Proof of Lemma~\ref{L2}:}
We notice that the triplet $\eta$ can be equivalently constructed
by choosing three points $(a,b,c)$ at random from $(N+1)$ points
on a circle and counting the number of points within each of the three
resulting segments.
This implies exchangeability. 

\qed

\begin{lem}
\label{L3}
Let $Y_i\in\mathbb{R}$, $i=1,2,\dots$ be i.i.d. random variables, a pair $(n,m)\in\mathbb{N}^2$
has an exchangeable distribution independent of $Y_i$, and
\be
X = \sum_{i=1}^{n} Y_i - \sum_{i=n+1}^{n+m} Y_i.
\ee
Then $X$ has a symmetric distribution.
\end{lem}

\noindent{\bf Proof of Lemma~\ref{L3}:}
Let $\Delta = n - m$ and $F(X\,|\,\Delta)$ denote the conditional 
distribution of $X$ given $\Delta$.
From the definition of $X$ it follows that 
\[F(X\,|\,\Delta = k) = F(-X\,|\,\Delta = -k).\]
Exchangeability of $(n,m)$ implies symmetry of $\Delta$
and we thus obtain
\begin{eqnarray}
F(X) &=& \sum _{k=-\infty}^{\infty} 
F(X\,|\,\Delta=k)\,{\sf P}(\Delta=k)\nonumber\\
&=& \sum _{k=0}^{\infty} 
\left[F(X\,|\,\Delta=k)+F(X\,|\,\Delta=-k)\right]
\,{\sf P}(\Delta=k)\nonumber\\
&=&\sum _{k=0}^{\infty} 
\left[F(X\,|\,\Delta=k)+F(-X\,|\,\Delta=k)\right]
\,{\sf P}(\Delta=k)\nonumber.
\end{eqnarray}
The sums of conditional distributions in brackets are symmetric, which 
completes the proof.
\qed

\vspace{.5cm}
\noindent
{\bf Proof of Theorem~\ref{T1}: Tokunaga self-similarity} 

We will show that 
$\displaystyle\lim_{N\to\infty}T_{ij} = 2^{j-i-1}$ for any pair $j>i$.
By Lemma~\ref{inv_prune}(b), $T_{ij} = T_{(i+k)\,(j+k)}$
and so it suffices to prove the statement for $i=1$, that is 
to show that $\displaystyle\lim_{N\to\infty}T_{1j} = 2^{j-2}$ for
any $j\ge 2$.
This will be done by induction.
Below we use the terminology introduced in Sect.~\ref{terms}.

{\it Induction base, $j=2$.}
Consider a basin of order 2, formed by two consecutive points from 
$\mathcal{T}_2$ (local minima of second order). 
We denote here their positions 
by $L$ and $R$, $L<R$.
This part of the proof will consider only local minima from this interval;
they will be referred to as ``points''.

The highest local minimum, or point $c=c_k^{(2)}\in\mathcal{S}_2$ forms 
a vertex parental to two branches of order 1 with Tokunaga indices $\{11\}$;
in addition, a random number of local minima corresponds to internal 
vertices parental to side-branches with Tokunaga indices $\{1j\}$, $j>1$.
The number $N^{(L,R)}_{12}$ of vertices of index $\{12\}$ within $(L,R)$ equals 
the number of side-branch points $X_k$ that are higher than their opposite minimum
of second order:
\[N^{(L,R)}_{12} = \#\{L<k<R\,:\,X_k>M^{(2)}_k\}.\]

For each side-branch vertex $X_k$ we necessarily have  
$X_k<X_c$ since $X_c$ is maximal among the local minima.
Recall that the local minima form a SHMC. 
Hence, for a randomly chosen side-branch $X_k$ we have
\[X_c-X_k = \sum_{i=1}^{\xi'} Y_i,\]
where $\xi'$ is a geometric rv such that
${\sf P}(\xi'=k)=2^{-k}$, and
$Y_i>0$ are i.i.d. random variables that correspond to the jumps between 
the local minima.
Clearly, the difference $X_c-M^{(2)}_k$ has the same distribution.
The random variables $(X_c-M^{(2)}_k)$ and $(X_c-X_k)$ are independent
and so ${\sf P}\left(X_k>M^{(2)}_k\right)=1/2$.
The expected number of side-branches with index $\{12\}$ within the 
interval $(L,R)$ is
\be
\label{rs}
{\sf E}\left[N^{(L,R)}_{12}\right] = 
{\sf E}\left[\sum_{k=1}^{\xi-1} \1_{(0,\infty)}\left(X_k-M^{(2)}_k\right)\right].
\ee
The summation above is taken over $(\xi-1)$ side-branch points within $(L,R)$;
and the random variables $\xi$ was described in Lemma~\ref{L4}.

We show next that the random variables  
$\1_{(0,\infty)}\left(X_k-M^{(2)}_k\right)$ are 
independent of $\xi$.
Suppose that there exist $\xi=N$ points within $(L,R)$.
A particular placement of $k$ and $c$ among these points is obtained by 
choosing two points at random and without replacement from $\{1,\dots,N\}$.
By Lemma~\ref{L2}, the conditional distribution of the numbers of 
points between $k$ and $c$ and between $c$ and the local minimum 
opposite to $X_k$ have an exchangeable distribution.  
Lemma~\ref{L3} implies that
${\sf P}\left(X_k>M^{(2)}_k\,|\,\xi=N\right)=1/2$.
Thus,
\be
{\sf E}\left[N^{(L,R)}_{12}\right]
={\sf E}[\xi-1]\,{\sf P}\left(X_k>M^{(2)}_k\right) = 2\times 1/2 =1.
\ee
The numbers $N_{12}^{(L,R)}$ are independent
for different basins of order 2 by Markov property of $X_t$.
The strong law of large numbers yields
\[T_{12} = \frac{N_{12}}{N_2}
\overset{\rm a.s.}{\longrightarrow}1=2^0
~{\rm as~}N\to\infty.\]

{\it Induction step.}
Suppose that the statement is proven for $j\ge 2$, that is we know that
for a randomly chosen local minima $X_k$
\[{\sf P}\left(X_k > M^{(j)}_k\right) = 2^{-(j-1)}\]
and $T_{1j}\overset{\rm a.s.}{\to}2^{j-2}$ as $N\to\infty$.
We will prove it now for $(j+1)$.
Consider a randomly chosen side-branch point $X_k$ of order $\{1i\}$, $i>j$.
By \eqref{sb}, $X_k < M^{(m)}_k$ for $1\le m \le j$ and thus 
necessarily $X_k < c_k^{(i+1)}$, $1\le i\le j$, since $c_k^{(i+1)}$ is a local
maximum of order-$i$ minima within the basin $(L,R)$ of 
order $(j+1)$ that contains $k$.
Repeating the argument of the induction base we find that
$X_k-M_k^{(i)}$ has a symmetric distribution for all $i\le j+1$
and that the probability of $(X_k>M_k^{(i)})$ is independent
of the number of local maxima of order $j$ within the basin $(L,R)$.
This gives, for a randomly chosen $X_k$,
\begin{eqnarray*}
{\sf P}\left(X_k > M^{(j+1)}_k\right) &=&
{\sf P}\left(X_k>M^{(j+1)}_k,X_k>M^{(j)}_k\right)\\
&=& {\sf P}\left(X_k > M^{(j+1)}_k \left| X_k > M^{(j)}\right.\right)\,
{\sf P}\left(X_k > M^{(j)}_k\right) \\
&=& 2^{-1}\times 2^{-(j-1)}=2^{-j}.
\end{eqnarray*}

By Lemma~\ref{L4}, the average number of order-2 basins within
a basin of order $(j+1)$ is $4^{j-1}$. 
Each such basin contains on average 2 points that correspond
to side branches with Tokunaga index $\{1\bullet\}$.
Hence, the average total number of side-branches with index $\{1\bullet\}$
within a basin of order $(j+1)$ is $2\times 4^{j-1}= 2^{2j-1}$.
Applying the Wald's lemma to the sum of indicators
$\1_{(0,\infty)}(X_k-M^{(j+1)}_k)$ over the random number of 
local minima of order $j$ within the basin $(L,R)$, we find
the average total number of side-branches of order $\{1(j+1)\}$: 
\[{\sf E}\left[N^{(L,R)}_{1(j+1)}\right] = 2^{-j}\times 2^{2j-1} = 2^{j-1}.\]
The strong law of large numbers yields
\[T_{1(j+1)} = \frac{N_{1(j+1)}}{N_{(j+1)}}
\overset{\rm a.s.}{\longrightarrow}2^{j-1},
~{\rm as~}N\to\infty.\]
\qed

{\bf Proof of Proposition~\ref{DSS}:}
Each transition step between the local minima of $X_k$ can be represented as 
$d_j$ of \eqref{twosum}
where $\{Y_i\}$ and $\{Z_i\}$ are independent random variables with density $f(x)$, 
and $\xi_+$ and $\xi_-$ are two independent geometric random variables 
with parameter $1/2$. 
The Wald's lemma readily implies that $c=2$.
This gives for the characteristic functions 
\[\widehat{K}_1(s)=2\,\widehat{K}(2s)=\Re\left[\widehat{f}(2s)\right].\]
On the other hand, taking the characteristic function of $d_j$
we obtain
\[\widehat{K}_1(s)=\left|\frac{\widehat{f}(s)}{2-\widehat{f}(s)}\right|^2,\]
which completes the proof.
\qed

\vspace{.5cm}
\noindent
{\bf Proof of Theorem~\ref{eH}:}
The Tokunaga and Horton self-similarity for a symmetric EHMC 
was proven in Theorem~\ref{T1}.
Here we show the violation of the Horton self-similarity
for an asymmetric EHMC.

Let $X^{(m)}_k$ denote the time series obtained by $m$-time repetitive
pruning of time series $X_k$.
Recall that there is one-to-one correspondence between the local maxima 
of $X^{(m)}$ and the branches of order $m$ in the level set tree 
$\textsc{level}(X)$ 
(see Sect.~\ref{terms}).
Hence, the Horton self-similarity is equivalent to the invariance of the 
proportion of local maxima with respect to pruning.
The proportion of local maxima in $X^{(m)}$ equals the probability 
$P^{(m)}_{\rm min}$ 
for a randomly chosen point to be a local maxima.  
The Markov property of $X^{(m)}$ --- Lemma~\ref{inv_prune}(c) --- implies that 
$P^{(m)}_{\rm min}=p^{(m)}(1-p^{(m)})$, where $p^{(m)}$ is
the probability for an upward jump in $X^{(m)}$. 

For an asymmetric EHMC let $A^{(m)}$ be the $m$-th iteration of $A$, 
as in \eqref{ag}, \eqref{dyn}. 
There, for $m\ge 1$, either 
$A^{(m)}<1$ in which case $A^{(m)}\to 0$ 
or 
$A^{(m)}>1$ in which case $A^{(m)}\to \infty$,
all as $m\to\infty$
(see Sect.~\ref{expj}, Eq.~\eqref{dyn} and Fig.~\ref{fig4}).
This corresponds to $p^{(m)}=1/(A^{(m)}+1)\to 1$ or $p^{(m)}\to0$, 
respectively, and leads to $P^{(m)}_{\rm min}\to0$.
This prohibits the Horton, and hence Tokunaga, self-similarity.
\qed

\section{Discussion}
\label{discussion}
This work establishes the Tokunaga and Horton self-similarity for
the level-set tree of a finite symmetric homogeneous Markov 
process with discrete time and continuous state space 
(Sect.~\ref{main}, Theorem~\ref{T1}).
We also suggest a definition of self-similarity for an infinite 
tree, using the construction of a forest of subtrees
attached to the floor line \cite{Pitman}; this allows us
to establish the Tokunaga and Horton self-similarity for
a regular Brownian motion (Sect.~\ref{main}, Corollary~\ref{Brown}).
This particular extension to infinite trees seems
natural for {\it tree representation of time series},
where concatenation of individual finite time series
corresponds to the ``horizontal'' growth of the corresponding
tree. 
Alternative definitions might be better 
suited though for other situations related, say, to the ``vertical'' 
growth of a tree from the leaves, like in a branching process.

A useful observation is the equivalence of smoothing the time 
series by removing its local maxima and pruning the
corresponding level-set tree  (Sect.~\ref{main}, Proposition~\ref{ts_prune}).
It allows one to switch naturally between the tree and 
time-series domains in studying various self-similarity 
properties.

As discussed in the introduction, the Tokunaga self-similarity 
for various finite-tree 
representations of a Brownian motion follow from 
(i) the results of Burd {\it et al.} \cite{BWW00} 
on the Tokunaga self-similarity for the critical binary 
Galton-Watson process and (ii) equivalence of a particular 
tree representation to this process.
We suggest here an alternative, direct approach to
establishing Tokunaga self-similarity in Markov processes.
Not only this approach does not refer to the Galton-Watson
property, it extends the Tokunaga self-similarity to a much 
broader class of trees. 
Indeed, as shown by Le Gall \cite{LeGall93} and Neveu and Pitman 
\cite{NP89} (see Theorem~\ref{Pit7_3}), the tree representation
of any non-exponential symmetric Markov chain is {\it not} 
Galton-Watson; it is still Tokunaga, however, by our Theorem~\ref{T1}. 

Peckham and Gupta \cite{PG99} have introduced the 
{\it generalized Horton laws}, which state the equality 
in distributions for the rescaled versions of suitable 
branch statistics $S_r$: 
$S_r \overset{d}{=} R_S^{r-k}\,S_k$, $R_S>0$.
These authors established the existence of the generalized
Horton laws in the Shreve's random model, that is for 
the Galton-Watson trees.
Accordingly, one would expect the generalized Horton laws to hold 
for the exponential symmetric Markov chains.
Veitzer and Gupta \cite{VG00} and Troutman \cite{Tr05} have 
studied the {\it random self-similar network (RSN) model} 
introduced in order to explain the variability of the limiting 
branching ratios in the empirical Horton laws.
They have demonstrated that the extended Horton laws hold for 
various branch statistics, including the average magnitudes $M_r$, 
in this model.
Furthermore, they established the weak Horton laws \eqref{NHL},
\eqref{MHL} and Tokunaga self-similarity for the RSN model.
Notably, the RSN model does not belong to the class of
Galton-Watson trees, yet it demonstrates the Tokunaga self-similarity, 
similarly to the non-exponential symmetric Markov chains considered here.

Tree representation of stochastic processes 
\cite{Pitman,NP89,Ald1,Ald2,Ald3,LeGall93,LeGall05} and 
real functions \cite{Arnold,EHZ03} is an intriguing topic that attracts 
attention of mathematicians and natural scientists.
A structurally simple yet flexible Tokunaga self-similarity,
which extends beyond the classical Galton-Watson space, 
may provide a useful insight into the structure of existing
data sets and models as well as suggest novel ways of modeling 
various natural phenomena.
For instance, the level set tree representation have been used recently in  
analysis of the statistical properties of fragment coverage in genome 
sequencing experiments \cite{EHP11,EHP10,Evans05}.
It seems that some of the methods and results obtained in this work
might prove useful for the gene studies.
In particular, it looks intriguing to test the self-similarity
of the gene-related trees and interpret it in the biological context.

Notably, the results of this paper, as well as that of 
Burd {\it et al.} \cite{BWW00}, refer only to a single point
$(a,c)=(1,2)$ in the two-dimensional space of Tokunaga parameters. 
The empirical and numerical studies, however, report a broad 
range of these parameters, roughly $1< a < 2$ and $1< c < 4$.
This motivates a search for more general Tokunaga models;
a potential broad family is suggested by our Conjecture~\ref{BH}.

The construction of the level set tree is a particular case of
the coagulation process; in the real function
context it describes the hierarchical structure of the 
embedded excursions of increasing lengths and heights.
Coagulation theory --- a well-established field with broad range
of practical applications to physics, biology, and social sciences
\cite{Bertoin,Wakeley,NBW06} --- is heavily based on the concepts
of symmetry and exchangeability \cite{Pitman,Bertoin}.
We find it noteworthy that the only property used to establish
the results in this paper is symmetry of a Markov chain.
It seems worthwhile to explore the concept of Tokunaga
self-similarity for a general coalescent process.

\vspace{0.5cm}
{\bf Acknowledgement.}
We are grateful to Ed Waymire and Don Turcotte for providing 
continuing inspiration to this study. 
We also thank Mickael Chekroun, Michael Ghil, 
Efi Foufoula-Georgiou, and Scott Peckham for their support
and interest to this work.
Comments of two anonymous reviewers helped us to 
significantly improve and expand an earlier version of 
this work. 
This study was supported by the NSF Awards DMS 0620838
and DMS 0934871.

\begin{figure}[p] 
\centering\includegraphics[width=0.8\textwidth]{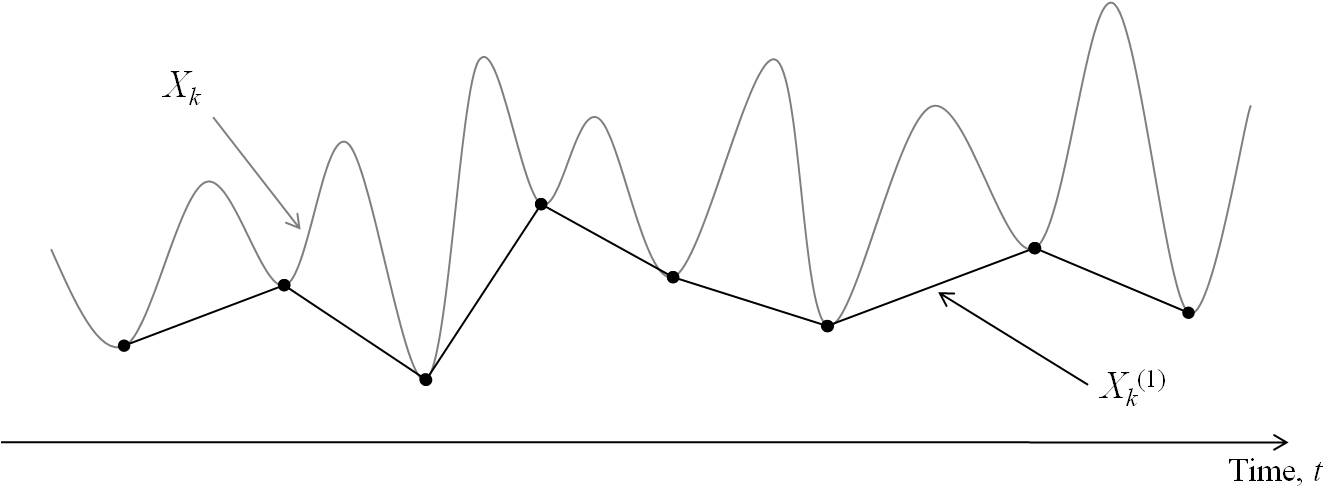}
\caption{Time series $X_k$ (light line) and the series $X_k^{(1)}$ 
of linearly connected local minima (black line and dots).}
\label{fig_ts_prune}
\end{figure}

\begin{figure}[p] 
\centering\includegraphics[width=0.8\textwidth]{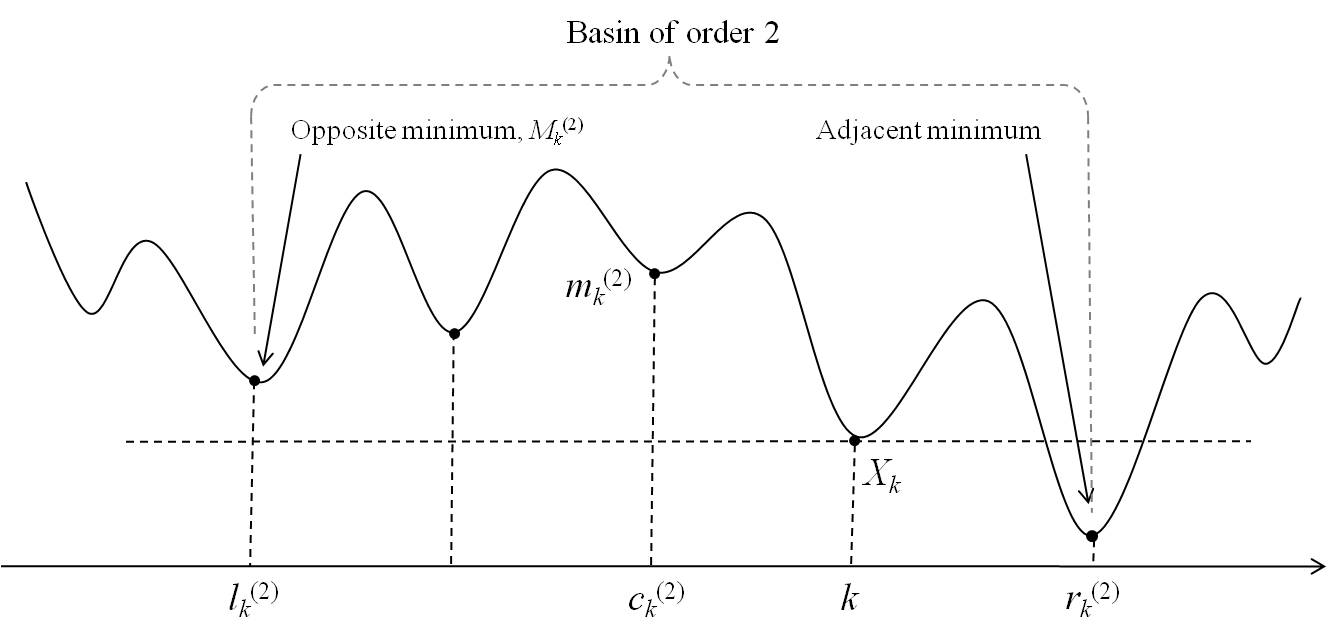}
\caption{Basin of order 2: an illustration.
The figure shows a basin of order 2 that consists of 5 local
minima.
The figure illustrates the taxonomy used in the paper; 
it shows the local maximum $m_k^{(2)}$ of the basin's local minima,
the opposite and adjacent minima of second order for a local 
minimum $X_k$,
as well as the corresponding points $l^{(2)}_k$, $c^{(2)}_k$, and 
$r^{(2)}_k$.
}
\label{fig1}
\end{figure}

\begin{figure}[h] 
\centering\includegraphics[width=0.7\textwidth]{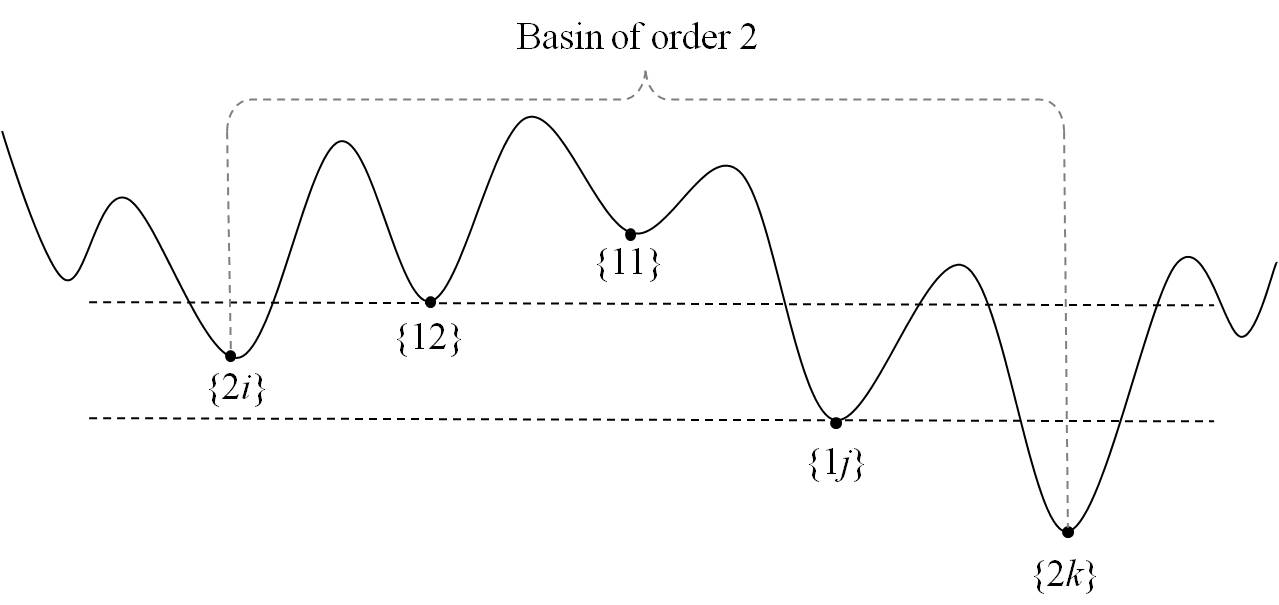}
\caption{Tokunaga indexing: an illustration. The figure shows the Tokunaga
indexing for the local minima of the second order basin shown in Fig.~\ref{fig1}.
The values of $i,j,k>2$ are determined by the large-scale structure of
the function $X_t$.}
\label{fig2}
\end{figure}

\begin{figure}[p] 
\centering\includegraphics[width=0.7\textwidth]{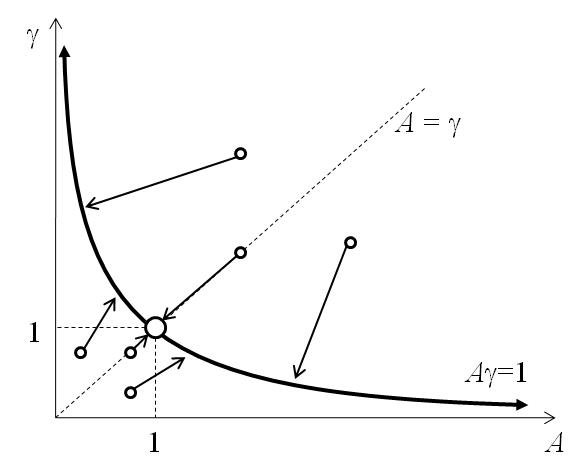}
\caption{Characterization of EHMCs in the space $(A,\gamma)$
of \eqref{ag} with iteration rules \eqref{dyn} that correspond to
the transition to the EHMC of local maxima. 
Each EHMC corresponds to a point on the plane $(A,\gamma)$.
The chain of local minima for any EHMC corresponds to
a point on the hyperbola $A\,\gamma=1$.
The point $(A=1,\gamma=1)$ is fixed. 
Any point from the lower branch ($A>1,\gamma<1$) moves
along the hyperbola toward $(\infty,0)$.
Any point from the upper branch ($A<1,\gamma>1$) moves
along the hyperbola toward $(0,\infty)$.
Arrows illustrate the point dynamics.
}
\label{fig4}
\end{figure}

\end{document}